\documentclass[11pt]{article}
\usepackage{amsfonts}
\usepackage{amssymb}
\begin{document}
\newtheorem{lemma}{Lemma}[section]
\newtheorem{theorem}{Theorem}[section]
\newtheorem{proposition}{Proposition}[section]
\newtheorem{corollary}{Corollary}[section]
\def\dim{\mathop{\rm dim}\nolimits}
\def\deg{\mathop{\rm deg}\nolimits}
\def\codim{\mathop{\rm codim}\nolimits}
\def\char{\mathop{\rm char}\nolimits}
\def\rank{\mathop{\rm rank}\nolimits}
\def\rn{\mathop{\rm rk}\nolimits}
\def\id{\mathop{\rm id}\nolimits}
\def\rnk{\mathop{\overline{\rm rk}}\nolimits}
\def\dg{\mathop{\rm dg}\nolimits}
\def\ind{\mathop{\rm ind}\nolimits}
\def\Hom{\mathop{\rm Hom}\nolimits}
\def\Hilb{\mathop{\rm Hilb}\nolimits}
\def\Im{\mathop{\rm Im}\nolimits}
\def\GL{\mathop{\rm GL}\nolimits}
\def\cA{\mathop{\cal A{}}\nolimits}
\def\cB{\mathop{\cal B{}}\nolimits}
\def\cD{\mathop{\cal D{}}\nolimits}
\def\cV{\mathop{\cal V{}}\nolimits}
\def\cQ{\mathop{\cal Q{}}\nolimits}
\def\cP{\mathop{\cal P{}}\nolimits}
\def\cU{\mathop{\cal U{}}\nolimits}
\def\cO{\mathop{\cal O{}}\nolimits}
\def\cG{\mathop{\cal G{}}\nolimits}
\def\cK{\mathop{\cal K{}}\nolimits}
\def\cT{\mathop{\cal T{}}\nolimits}
\def\cZ{\mathop{\cal Z{}}\nolimits}
\def\cN{\mathop{\cal N{}}\nolimits}
\def\cM{\mathop{\cal M{}}\nolimits}
\def\cH{\mathop{\cal H{}}\nolimits}
\def\cX{\mathop{\cal X{}}\nolimits}
\def\cY{\mathop{\cal Y{}}\nolimits}
\def\cW{\mathop{\cal W{}}\nolimits}
\def\cI{\mathop{\cal I{}}\nolimits}
\def\cR{\mathop{\cal R{}}\nolimits}
\def\cC{\mathop{\cal C{}}\nolimits}
\def\cE{\mathop{\cal E{}}\nolimits}
\def\cF{\mathop{\cal F{}}\nolimits}
\def\cS{\mathop{\cal S{}}\nolimits}
\def\cSN{\mathop{\cal SN{}}\nolimits}
\def\cSE{\mathop{\cal SE{}}\nolimits}
\def\cSD{\mathop{\cal SD{}}\nolimits}
\def\cSC{\mathop{\cal SC{}}\nolimits}
\def\cSU{\mathop{\cal SU{}}\nolimits}
\def\pf{\mbox{\bf Proof \ }}
\def\pfA{\mbox{\bf Proof A \  }}
\def\pfB{\mbox{\bf Proof B \ }}
\newcommand{\N}{\mbox{$\mathbb N$}}
\newcommand{\C}{\mbox{$\mathbb C$}}
\newcommand{\R}{\mbox{$\mathbb P$}}
\newcommand{\Z}{\mbox{$\mathbb Z$}}
\newcommand{\Q}{\mbox{$\mathbb Q$}}
\newcommand{\n}{\mbox{${\displaystyle \frac{n}{2} }$}}
\newcommand{\m}{\mbox{${\displaystyle \Big[ \frac{n}{2} \Big] }$}}
\newcommand{\mm}{\mbox{$ \Big[ \frac{n}{2} \Big] $}}
\hyphenation{as-so-cia-ted con-si-de-ring sti-mu-la-ting}
\title{On the maximum nilpotent orbit intersecting a centralizer in
$M(n,K)$}
\author{{Roberta Basili}
}
\date{}
\maketitle \noindent
\begin{abstract}   Recently
several authors have published results about the problem of
finding which pairs of partitions of $n$ correspond to pairs of
commuting $n\times n$ nilpotent matrices. We describe a maximal
nilpotent subalgebra $\cSN _B$ of the centralizer of a given
nilpotent $n\times n\, $ matrix with Jordan form; we show that
$\cSN _B$ intersects the maximum nilpotent orbit intersecting the
minimal affine subspace defined by the condition that some
coordinates are $0$ of which $\cSN _B$ is a subvariety. Then we
prove that the maximum partition which forms with a given partition
a pair with the previous property can be found by a
simple algorithm which was conjectured by Polona Oblak.
\vspace{1mm}
\vspace{1mm} \newline {\em Keywords:} commuting
matrices, nilpotent orbits, maximal orbit.
\end{abstract}
\section{Introduction} We will use the following notations: $
M (n,K) $ is  the set of the $n\times n$ matrices over an infinite
field $K$, $\, \mbox{\rm GL}\, (n,K) $ is  the set of the $n\times
n$ nonsingular matrices over $K$, $\, N (n,K)$ is the set of the
$n\times n$  nilpotent matrices over $K$, $\, J\, \in N (n,K)$ is
a matrix with Jordan canonical form, $\, \mu _1\geq \mu _2\geq
\cdots \geq \mu _t $ are  the orders of the Jordan blocks of $J$,
$\, B=(\mu _1,\ldots ,\mu _t)$ is the partition of $n$ associated
to $J$ and the orbit of $J$ under the action of GL$(n,K)$.\newline
Let $ J' $ be another nilpotent matrix, let $\mu '_1\geq \cdots
\geq \mu '_{t'}$ be the orders of its Jordan block and let $
B'=(\mu '_1,\ldots ,\mu '_{t'})$. Then $B=B'$ if and only if
$\rank J^m=\rank (J')^m$ for all $m\in \N $. It is said that $
B<B'$ if $\, \rank J^m\leq \rank (J')^m\, $ for all $ m\in \N$ and
there exists $ m\in \N$ such that $\, \rank J^m\, <\, \rank
(J')^m$. The following claim is due to Hesselink (\cite{He},
1976): $B<B'\,  $ iff $  B$, as an orbit, is contained in the
closure of the orbit $ B'$. \newline For $i\in \N $ we set $\mu
_i=0$ for $i>t$, $\mu'_i=0$ for $i>t'$; then we have that
$$ B\leq B' \qquad \Longleftrightarrow \qquad {\displaystyle \sum_{i=1}^l\mu _i \leq
\sum_{i=1}^l\mu '_i}\quad \mbox{\rm for all } \ l\in \N
$$ where equality holds in the first relation iff it holds in the
second one for all $l\in \N $. \vspace{2mm}
\newline {\bf Examples}
 $ (6,4,3)<(6,5,2)<(6,6,1)$, $\ (5,3,2,1)<(6,3,1,1)<(6,4,1)
\ ; $ $(6,5,2)$ and $(7,3,3)$ cannot be compared, as $(6,5,4,3)$
and $(6,6,2,2,2)$. \vspace{2mm}\newline Let $ {\cal C }_B$ be the
centralizer of $J$ and let $ \cN _B$ be the subset of $\cC _B$ of
all the nilpotent matrices. We recall the following result, whose
proof is a consequence of Wedderburn's Theorems.
\begin{lemma}\label{algebra} If $\cU $ is a finite dimensional
algebra over an infinite field $K$ then the subvariety $\cN (\cU)$
of the nilpotent elements of $\cU $ is irreducible.\end{lemma} For
$m\in \N $ the subvariety of $\cN _B$ of all $X$ such that $\rank
X^m$ is the maximum possible is open and by lemma \ref{algebra}
(see also for example Lemma 2.3 of \cite{Bas}) $\cN _B$ is
irreducible, then the intersection of these open subsets for $m\in
\N $ is not empty.
 Hence there is a maximum partition for the elements of
$\cN _B$ and the subset of the elements which have this partition
is open (dense) in $\cN _B$.
 Then we can define a map $Q$ in the set of
the orbits of $n\times n$ nilpotent matrices (or partitions of
$n$) which associates to any orbit $B$ the maximum nilpotent orbit
which intersects $\cN _B$.\newline  For $\ s\in \N -\{ 0\} \ $ let
$q$ and $r$ be the quotient and the remainder of the division of
$n$ by $s$; then if $J$ is the $n\times n$ Jordan block we have
that $J^s$ has $r$ Jordan blocks of order $q+1$ and $s-r$ Jordan
blocks of order $q$. Hence the partition $B$ is almost rectangular
(that is $\mu _1-\mu _t\leq 1$) iff $J$ is conjugated to a power
of the $n\times n$ Jordan block. Hence if $B$ is almost
rectangular we have $Q(B)=(n)$. The converse of this claim is also
true; it is a consequence of the next Proposition which we are
going to explain.
\newline There exist $p\in \N $ and almost rectangular
partitions $B_1,\ldots ,B_p$ of numbers less or equal than $n$
such that $B=(B_1,\ldots ,B_p)$; we denote by $r_B$ the minimum of
all $p$ with this property (sometimes it can be obtained with
different choices of $B_1,\ldots ,B_{r_B}$). \vspace{2mm}\newline
{\bf Examples } If $B=(5,4,3,1,1)$ we have $r_B=3$, which is
obtained with $B_1=(5)$, $B_2=(4,3)$, $B_3=(1,1)$ or with
$B_1=(5,4)$, $B_2=(3)$, $B_3=(1,1)$. If $B=(9,7,5,1)$ we have
$r_B=4$. \vspace{2mm}
\newline Let $s_B$ be the maximum value of $l$ for which there
exists
 a subset $\{ i_1,\ldots , i_l\} $ of $ \{ 1,\ldots ,t\} $ such that
 $i_1< \ldots <i_l$ and $\mu _{i_1} -\mu _{i_l} \leq 1$. The
 following Propositions were proved in {\cite{Bas} (Propositions 2.4
 and 3.5).
\begin{proposition} \label{R1} There exists a non-empty open subset of $\cN
_B$ such that if $A$ belongs to it we have that $\mbox{\rm rank}\
A=n-r_B$ (that is $A$ has $r_B$ Jordan blocks).
\end{proposition}
\begin{proposition}\label{R2} For $A\in \cN _B$ and $m\in \N
$ we have  that $$\mbox{\rm rank }\  (A^{s_B})^m\leq \mbox{\rm
rank } \ J^m\ .$$\end{proposition} Let $B=(B_1,\ldots ,B_{r_B})$
where $B_i$ is almost rectangular for $i=1,\ldots , r_B$; let
$n_i$ be the sum of the numbers of $B_i$ and let $\widetilde B
=(n_1,\ldots , n_{r_B})$. The following result is a consequence of
Proposition \ref{R2} (see Theorem 1.11 of  \cite{BasI}).
\begin{proposition}\label{RI} If $\, s_B=|B_i| $ for $ \
i=1,\ldots ,r_B\ $ then  $\ Q(B)=\widetilde B$.\end{proposition}
{\bf Example} \hspace{2mm} If $B=(5,4,4,2,2,1)$ we have
$\widetilde B=(13,5)$ and $Q(B)=\widetilde B$.
\begin{corollary}\label{corR} We have that $\ Q(B)=B\ $ iff $\ r_B=t\ $, that is
 $ \  \mu _{i}-\mu _{i+1}>1\ $  for
$ \ i=1,\ldots ,t-1$. \end{corollary}
 Let $\{ q_1,q_2,\ldots ,q_u\} $ be the
ordered subset of $\{ 0,\ldots ,t\} \ $  such that $\ q_u=t\ $ and
$$\mu _1=\mu _{q_1}\neq \mu _{q_{1}+1}=\mu _{q_{2}}\neq \mu _{q_2+1}=\cdots \neq
\mu _{q_{u-1}+1}= \mu _{q_u}\ $$  (for example if $B=(6,6,6, 6,
5,2, 2, 1)$ we have that   $q_1=4$, $q_2=5$, $q_3=7$, $q_4=8$). If
we set $q_0=0$ then $J$ has  $ q_i-q_{i-1}$ Jordan blocks of order
$ \mu _{q_i}$ for $i=1,\ldots ,u$.
 We
will write the partition $(\mu _1,\ldots ,\mu _t)$ also as $$ (\mu
_{q_1}^{q_1},\mu _{q_2}^{q_2-q_1}\ldots ,\mu _{q_{u}
}^{q_{u}-q_{u-1}}) \ .$$ We consider the subset of  $ \{ 1,\ldots
,u\} \times \{ 0,1\} $ of all $(i, \epsilon )$ such that $\mu
_{q_i}-\mu _{q_{i+\epsilon }}\leq 1$ (that is such that $\epsilon
=1$ iff $i<u$ and $\mu _{q_{i+1}}=\mu _{q_i}+1$); then we consider
the map from this subset to $\N $ defined by
$$(i, \epsilon ) \quad \longmapsto \quad 2q_{i-1}+\mu
_{q_i}(q_i-q_{i-1})+\epsilon \, \mu _{q_{i+1}}(q_{i+1}-q_{i})\ .$$
We denote by $\omega _1$ the maximum of the image of this map and
by $(\tilde i, \tilde {\epsilon })$ a preimage of $\omega _1$.
Polona Oblak in \cite{Obl} (2007) proved the following result.
\begin{theorem} \label{O} The maximum index of nilpotency for an element of $\cN
_B$ (that is the first number of the partition $Q(B)$) is $\
\omega _1$.
\end{theorem}{\bf Example} \ If $B=(5^2,4,3^4,2,1)$ we can only have $\tilde i=2$, $\tilde {\epsilon
}=1$ and the first number of $Q(B)$ is $\omega_1=2\times 2 +
4+3\times 4=20$.\vspace{2mm} \newline We will write the canonical
basis $\Delta _B $ of $K^n$ in the following way:
$$\Delta _B=\{ \ v_{\mu _{q_i},j}^{\mu _{q_i}}, v_{\mu _{q_i},j}^{\mu _{q_i}-1},\ldots
, v_{\mu _{q_i},j}^{1},\ i=1,\ldots ,u,\ j=q_i-q_{i-1},\ldots ,1
\} \ .$$ For example if $B=(5,3,3,2,1)=(5^1,3^2,2^1,1^1)$ we will
write
$$\Delta _B=\{ \ \underbrace{v_{5,1}^5,v_{5,1}^4,\ldots
,v_{5,1}^1},\underbrace{v_{3,2}^3,v_{3,2}^2, v_{3,2}^1},
\underbrace{v_{3,1}^3, v_{3,1}^2, v_{3,1}^1},\underbrace{
v_{2,1}^2, v_{2,1}^1},\underbrace{ v_{1,1}^1}\ \} \
.\vspace{2mm}$$ Let $\Delta _B^{\circ }$ be the union of the
following subsets:
$$\Delta _B^{\circ ,1}=\{ v_{\mu _{q_i},j}^1\ |\ j=q_i-q_{i-1},\ldots ,1\ , \ i=1,\ldots
,\tilde i+\tilde{\epsilon}\} \ , $$
$$\Delta _B^{\circ ,2}=\{
v_{\mu _{q_i},j}^l \ | \ j=q_i-q_{i-1}, \ldots ,1 \ ,\ l=1,\ldots
,\mu _{q_i} \ ,\ i=\tilde i,\tilde i+\tilde {\epsilon }\} \ ,$$
$$\Delta _B^{\circ ,3}=\{ v_{\mu
_{q_i},j}^{\mu _{q_i}}\ |\ j=q_i-q_{i-1},\ldots ,1\ , \ i=1,\ldots
,\tilde i+\tilde{\epsilon} \} \ . $$ Then $\big| \Delta _B^{\circ
}\big| =\omega _1$.  Let $\widehat B$ be the partition obtained
from $B$ by cancelling the powers $\mu _{q_{i}}^{q_{i}-q_{i-1}}$
for $i=\tilde i,\tilde i+\tilde{\epsilon}$ and decreasing by $2$
the numbers $\ \mu _{q_i}$ for $i=1,\ldots ,\tilde i-1$, that is:
$$\widehat B= \big( (\mu _{q_1}-2)^{q_1}, \ldots ,(\mu
_{q_{\tilde i-1}}-2)^{q_{\tilde i-1}-q_{\tilde i-2}}, \mu
_{q_{\tilde i+\tilde{\epsilon}+1}}^{q_{\tilde
i+\tilde{\epsilon}+1}-q_{\tilde i+\tilde{\epsilon}}}, \ldots , \mu
_{q_u}^{q_u-q_{u-1}}\big)  \ .$$  Let $Q(B)=(\omega_1,\ldots
,\omega_z)$; let $Q(\widehat {B})=(\widehat {\omega } _1, \ldots
,\widehat{\omega }_{\hat z})$. The main aim of this paper is to
prove the conjecture of Polona Oblak which is expressed by the
following Theorem.
\begin{theorem}\label{eend} The maximum partition $(\omega_1,\omega_2,\ldots ,\omega_z)$ which is associated to
 elements of $ {\cN }_B$ is  $(\omega_1, \widehat{\omega } _1, \widehat{\omega}_2,\ldots , \widehat
{\omega}_{\hat z})$, that is $Q(B)=(\omega_1, Q(\widehat{B}))$.
\end{theorem} This conjecture was communicated by its author to some of the participants
of the meeting "Fifth Linear Algebra Workshop" (Kranjska Gora, May
27 - June 5, 2008). Recently A. Iarrobino and L. Khatami wrote a
paper on the inequality $(\omega_1,\omega_2,\ldots ,\omega_z)\geq
(\omega_1, \widehat{\omega } _1, \widehat{\omega}_2,\ldots ,
\widehat {\omega}_{\hat z})$ (see \cite{KhI}; the first version of
it was circulated in October 2011). In that paper the conjecture is
expressed by using the digraph associated to the nilpotent
subalgebra $\cSN _B$ which will be defined in the next section
(see also \cite{Obl}); by Theorem \ref{O} that form is equivalent
to Theorem \ref{eend}. Recently L. Khatami obtained the following results: in \cite{Kha} she proved the uniqueness of the result of the
algorithms which follow from Theorem \ref{eend} (since there can
be different choices of $(\tilde i ,\tilde {\epsilon })$), which
we will show to be $Q(B)$; in \cite{Khat} she gave a formula for
$\omega_z$.
\section{The subalgebras ${\cE }_B$, ${\cN }_B$,
$\cSE _B$ and $\cSN _B$} We will consider any $n\times n$ matrix
$X$ as a block matrix $(X_{h,k})$, where $X_{h,k}$ is a $\mu
_h\times \mu _k$ matrix and $h,k=1,\ldots ,t$. Let $ \cU _B$ be
the subalgebra of $M(n,K)$ of all $X$ such that for $1\leq k\leq
h\leq t$ the blocks $X_{h,k}$ and $X_{k,h}$ are upper triangular,
that is have the following form:\vspace{2mm}
$$X_{h,k}=\pmatrix{0 &\ldots & 0 & x_{h,k}^{1,1} & x_{h,k}^{2,1} & \ldots &
x_{h,k}^{\mu _h,1} \cr \vdots &  &  & 0 & x_{h,k}^{1,2} & \ddots &
x_{h,k}^{\mu _h-1,2} \cr \vdots &  &  &  & \ddots & \ddots &
\vdots \cr 0 & \ldots & \ldots & \ldots & \ldots & 0 &
x_{h,k}^{1,\mu _h} \cr},\vspace{2mm}$$
$$X_{k,h}=\pmatrix{x_{k,h}^{1,1} & x_{k,h}^{2,1} & \ldots & x_{k,h}^{   \mu _h,1} \cr
0 & x_{k,h}^{1,2} & \ddots & x_{k,h}^{\mu _h-1,2} \cr \vdots &
\ddots & \ddots & \vdots \cr \vdots &  & 0 & x_{k,h}^{1,\mu _h}
\cr \vdots & & & 0 \cr \vdots &  &  & \vdots \cr 0 & \ldots &
\ldots & 0 \cr}\vspace{2mm}$$ where for $\mu _h=\mu _k$ we omit
the first $\mu _k-\mu _h$ columns and the last $\mu _k-\mu _h$
rows respectively. \newline For $X\in \cU _B $, $\ i,j \in \{
1,\ldots ,u\} $ and $l\in \{ 1,\ldots ,\mbox{\rm min }\{ \mu
_{q_{i }}, \mu _{q_{j }}\} \} $ let
$$X (i,j ,l)=(x_{h,k}^{1,l})\ ,\qquad  q_{i
-1}+1\leq h\leq  q_{i},\ q_{j -1}+1\leq k\leq  q_{j}.$$ We set
$X(i,l)=X(i,i,l)$. Lemma \ref{algebra} for the algebra $\cU _B$
becomes more precise, as follows.
\begin{lemma} \label{1}
For $X\in \cU _B $ we have that: \begin{itemize}  \item[a)] there
exists $G\in $ GL $(n,K)$ such that $G^{-1}X G\,\in \,\cU _B$ and $ (G^{-1}X G)\; (i ,l) $ is
lower (upper) triangular  for $\, i =1,\ldots ,u\, $ and $\, l=
1,\ldots ,\mu _{q_{i }}\, $; \item[b)]
 $X$ is nilpotent iff $\, X( i , l)\, $ is nilpotent for $\, i =1,\ldots ,u\, $ and $\, l= 1,\ldots
,\mu _{q_{i }} \, $.\end{itemize}
\end{lemma}
\pf  We can construct a semisimple subalgebra of $\cU _B $ whose
direct sum with the Jacobson radical of $\cU _B$ is $\cU _B $; the
construction is as follows. For $l=1,\ldots ,\mu _{q_1}$ let
$$U^l=\langle v_{\mu _{q_i},j}^l \ : \ i=1,\ldots ,u,\
j=q_i-q_{i-1},\ldots ,1,\  \mu _{q_i}\geq l\rangle ;$$ then
$K^n={\displaystyle \bigoplus _{l=1}^{\mu _{q_1}}} U^l$ and
$X(U^l)\subseteq {\displaystyle \bigoplus _{i=l}^{\mu
_{q_1}}}U^i$. For $v\in K^n$ let $v={\displaystyle \sum
_{l=1}^{\mu _{q_1}}}v^{(l)}$ where $v^{(l)}\in U^l$ and let
$L_{X,l}\, :\, U^l\to U^l$ be defined by $L_{X,l}(v)=X(v)^{(l)}$.
Then $X$ is nilpotent iff $L_{X,l}$ is nilpotent for $l=1,\ldots
,\mu _{q_1}$. For $l=1,\ldots ,\mu _{q_1}$ let $i _l\in \{
1,\ldots ,u \} $ be such that $l\leq \mu _{q_{i_l}}$ and $\mu
_{q_{i_l+1}}<l$ if $i_l\neq u$. Then the matrix of $L_{X,l}$ with
respect to the basis $\{ v_{\mu _{q_i},j}^l\ \  : \ i=1,\ldots
,u,\ \mu _{q_i}\geq l\} $ is the lower triangular block matrix
$\big( X (i,j ,l)\big) $, $i ,j =1,\ldots ,i _l$, which is
nilpotent iff $X (i, l)$ is nilpotent for $i =1,\ldots ,i _l$. For
$v\in K^n$ let $v={\displaystyle \sum _{i = 1}^{u } v_{(i)}}$
where $v_{(i )}\in \langle v_{\mu_{q_i},j}^l \ :\
j=q_i-q_{i-1},\ldots ,1,\   l=1,\ldots ,u_{q_{i}}\rangle $. For $i
=1,\ldots ,u $ and $l=1,\ldots ,u_{q_{i }}$ let $U^l_{i}=\langle
v_{\mu_{q_i},j}^l\ :\ j=q_i-q_{i-1},\ldots ,1\rangle $ and let $
L_{X,i ,l}\ :\ U^l_{i}\to U^l_{i} $ be defined by $L_{X,i
,l}(v)=L_{X,l}(v)_{(i)}$. Then $ X(i ,l )$ is the matrix of $L
_{X,i,l}$ with respect to the basis $\{ v_{\mu_{q_i},j}^l\ :\
j=q_i-q_{i-1},\ldots ,1\} $ . We can substitute this basis with
another basis of the same subspace such that $X(i,l)$ is upper
triangular, for $i=1,\ldots ,u$ and $j=q_i-q_{i-1},\ldots
,1$.\hspace{4mm} $\square$ \vspace{2mm}\newline
 We
will denote by $ \cSU   _B $ the subspace of $\cU _B $ of all $X$
such that $X(i ,l)$ is lower triangular for $i =1,\ldots ,u$ and
$l=1,\ldots ,\mu _{q_{i }}$. Moreover we will denote by $\cE _{B
}$ the subset of $\cU _B $ of all the nilpotent matrices and we
will set
$$ \cSE _{B }=\cSU   _B\; \cap \; \cE _{B }\ .$$
\begin{lemma} \label{2}
 The centralizer $\cC _B$ of $J$ has the following properties: \begin{itemize} \item[i)] it is the
subspace of  $\cU _B$ of all $X$ such that
$$x_{h,k}^{l,1}=x_{h,k}^{l,2}=\cdots =x_{h,k}^{l,\mu _h+1-l}\ ,
\quad
 \ x_{k,h}^{l,1}=x_{k,h}^{l,2}=\cdots = x_{k,h}^{l,\mu _h+1-l}$$
for $\ 1\leq k\leq h\leq t\ $ and $\ l=1,\ldots ,\mu _h$;
\item[ii)] if $X\in \cC _B$ we can choose $G$ with the property
expressed in a) of lemma \ref{1} and such that
$GJ=JG$.\end{itemize}
\end{lemma}
\pf For i) see \cite{Ait} or Lemma 3.2 of \cite{Bas'}. Using the
notations of the proof of lemma \ref{1}, for $i=1,\ldots ,u$ let
${\big (} c^{(i)}_{h,k}{\big )} $, $h,k=q_i-q_{i-1},\ldots ,1$ be
a $q_i-q_{i-1}$ matrix over $K$ such that the vectors
$$w_{\mu _{q_i},j}^1={\displaystyle \sum
_{k=q_i-q_{i-1}}^1c^{(i)}_{j,k}v_{\mu _{q_i},k}^1}$$  form a basis
with respect to which $L_{X,i,1}$ is upper triangular. If we set
$$w_{\mu _{q_i},j}^l={\displaystyle \sum
_{k=q_i-q_{i-1}}^1c^{(i)}_{j,k}v_{\mu _{q_i},k}^l}$$ for
$l=1,\ldots \mu _{q_i}$ and for $i=1,\ldots ,u$ we get the basis
required by ii).\hspace{4mm} $\square $\vspace{2mm}\newline We can
shortly say that
 $X\in \cC _B$ iff its blocks are upper triangular
Toeplitz matrices.
 By lemma \ref{2} if $\, A\in \cC _B\, $ then $\, A
(i ,l)= A(i ,l') \, $ for $\, i \in \{ 1,\ldots ,u\} \, $ and $\,
l,l'\in \{ 1,\ldots ,\mu _{q_{i }} \} \, $; we denote this matrix
by $ \, A(i )\, $. \newline We will denote by $ \cSC _B$ the
subspace of all $A\in \cC _{B } $ such that $A(i )$ is lower
triangular for $i =1,\ldots ,u$. Moreover we will set
$$ \cSN _B=\cSC _B\; \cap \; \cN _B\
.$$ {\bf Example 1} If $B=(3,3,3,2)$ we have that $A\in \cN _B$
iff there exists a set  $\{ a_{h,k}^l\in K\  |\ (h,k,l)\in \{
1,2,3,4\} ^2 \times \{ 1,2,3\}\} $ such that $A$  is the
matrix:\vspace{2mm}
$${ {\pmatrix{ a_{11}^1 & a_{11}^2 & a_{11}^3 & | &  a_{12}^1 & a_{12}^2
& a_{12}^3 & | &  a_{13}^1 & a_{13}^2 & a_{13}^3 & | & a_{14}^1 &
a_{14}^2\cr  & a_{11}^1 & a_{11}^2 & | &   & a_{12}^1 & a_{12}^2 &
| &
  & a_{13}^1 & a_{13}^2 & | &  & a_{14}^1 \cr & & a_{11}^1 & | & & & a_{12}^1 & | & &  & a_{13}^1 & | & &
\cr - & -& -& -& -& -& -& -& -& -& - & -& -& - \cr a_{21}^1 &
a_{21}^2 & a_{21}^3 & | &  a_{22}^1 & a_{22}^2 & a_{22}^3 & | &
a_{23}^1 & a_{23}^2 & a_{23}^3 & | & a_{24}^1 & a_{24}^2\cr &
a_{21}^1 & a_{21}^2 & | & & a_{22}^1 & a_{22}^2 & | &
 &  a_{23}^1 & a_{23}^2 & | &  & a_{24}^1 \cr & & a_{21}^1 & | & & & a_{22}^1 & | & & & a_{23}^1 & | &  &
 \cr - & -& -& -& -& -& -& -& -&
-& - & -& -& - \cr a_{31}^1 & a_{31}^2 & a_{31}^3 & | & a_{32}^1 &
a_{32}^2 & a_{32}^3 & | & a_{33}^1 & a_{33}^2 & a_{33}^3 & | &
a_{34}^1 & a_{34}^2\cr  & a_{31}^1 & a_{31}^2 & | &  & a_{32}^1 &
a_{32}^2 & | &
 &  a_{33}^1 & a_{33}^2 & | &  & a_{34}^1 \cr & & a_{31}^1 & | & & & a_{32}^1 & | & & & a_{33}^1 & | & &
\cr - & -& -& -& -& -& -& -& -& -& - & -& -& - \cr & a_{41}^1 &
a_{41}^2 & | &  & a_{42}^1 & a_{42}^2 & | & & a_{43}^1 & a_{43}^2
& | & 0 & a_{44}^2\cr & &   a_{41}^1 & | & & & a_{42}^1 & | &  & &
a_{43}^1 & | & & 0 \cr
 }}}$$ and $$\pmatrix{a_{11}^1 & a_{12}^1 & a_{13}^1 \cr a_{21}^1 & a_{22}^1
& a_{23}^1 \cr a_{31}^1 & a_{32}^1 & a_{33}^1\cr} \in N(3,K)\ .$$
For each $A\in \cN_B$ there exists $\, G\in \cC _B\, $ and a set
$\{ \bar a_{h,k}^l\in K\  |\ (h,k,l)\in \{ 1,2,3,4\} ^2 \times \{
1,2,3\} \} $
 such
that det $G\neq 0$ and $ G^{-1}AG$ is the following element of
$\cSN _B$:\vspace{2mm}
 $${ {\pmatrix{ 0 & \bar a_{11}^2 & \bar a_{11}^3 & | & 0 &
\bar a_{12}^2 & \bar a_{12}^3 & | & 0 & \bar a_{13}^2 & \bar
a_{13}^3& | & \bar a_{14}^1 & \bar a_{14}^2\cr  & 0 & \bar
a_{11}^2 & | & & 0 & \bar a_{12}^2 & | &
 &  0 & \bar a_{13}^2 & | &  & \bar a_{14}^1 \cr & & 0 & | & &  & 0 & | & &  & 0 & |  & & \cr
- & -& -& -& -& -& -& -& -& -& - & -& -& - \cr \bar a_{21}^1 &
\bar a_{21}^2 & \bar a_{21}^3 &  | & 0 & \bar a_{22}^2 & \bar
a_{22}^3 & | & 0 & \bar a_{23}^2 & \bar a_{23}^3 & | & \bar
a_{24}^1 & \bar a_{24}^2\cr  & \bar a_{21}^1 & \bar a_{21}^2 & | &
& 0 & \bar a_{22}^2 & | &
 &  0 & \bar a_{23}^2 & | &  & \bar a_{24}^1 \cr & & \bar a_{21}^1 & | & & & 0  & | & & & 0 & | & &
 \cr - & -& -& -& -& -& -& -& -&
-& - & -& -& - \cr \bar a_{31}^1 & \bar a_{31}^2 & \bar a_{31}^3 &
| & \bar a_{32}^1 & \bar a_{32}^2 & \bar a_{32}^3 & | & 0 & \bar
a_{33}^2 & \bar a_{33}^3& | & \bar a_{34}^1 & \bar a_{34}^2\cr &
\bar a_{31}^1 & \bar a_{31}^2 & | & & \bar a_{32}^1 & \bar
a_{32}^2 & | &
 &  0 & \bar a_{33}^2 & | &  & \bar a_{34}^1 \cr & & \bar a_{31}^1  & | & &  & \bar a_{32}^1  & | & &  & 0  & |  &  &
 \cr
 - & -& -& -& -& -& -& -& -&
-& - & -& -& - \cr
 & \bar a_{41}^1 & \bar a_{41}^2 & | &  & \bar a_{42}^1 & \bar a_{42}^2 & | & & \bar a_{43}^1
& \bar a_{43}^2 & | & 0 & a_{44}^2 \cr & &   \bar a_{41}^1 & | & &
& \bar a_{42}^1 & | & & & \bar a_{43}^1 & | &  & 0  \cr
 }}}\vspace{2mm}$$
By lemma \ref{1} we get the following result.
\begin{corollary}\label{fromlemma2} The subvariety $\cSE
 _B$ ($\cSN _B$) has non-empty intersection with the orbit of any
element of $\cE _B$ ($\cN _B$).\end{corollary}
\section{Upper
triangular form of $ \cSE _B$, $ \cSN _B$}  We will denote by $<$
the linear order of $\Delta _B$ given by
  $v_{\mu _{q_i},j}^l \, < \,
v_{\mu_{q_{i'}},j'}^{l'} \ $  iff one of the following conditions
holds: \begin{itemize} \item[c$_1$)] $\ i<i'\ $ (that is $\ \mu
_{q_i}>\mu _{q_{i'}} $);\item[c$_2$)] $\ i=i'\ $ and $\ j>j'$;
\item[c$_3$)] $ \ i=i'\ $, $\ j=j'\ $ and $\ l>l'$. \end{itemize}
We observe that for $i\in \{ 1,\ldots ,u\} $ an element of $B$ is
equal to $\mu _{q_i}$ iff there exists $j\in \{ q_i-q_{i-1},\ldots
,1 \} $ such that it is the $ (q_i-j+1)-$th element of $B$;  if we
represent an element $X$ of $M(n,K)$ in the block form $
X=(X_{h,k})$, $\ h,k=1,\ldots ,t $ we have that the entry of $X\,
v_{\mu _{q_{i'}},j'}^{l'}$ with respect to $v_{\mu _{q_i},j}^l$ is
the entry of the matrix $X_{q_i-j+1, q_{i'}-j'+1}$ which, in this
matrix, has indices $\big( \mu _{q_i}-l+1, \mu _{q_{i'}}-l'+1\big)
$. The next lemma explains which square blocks of the matrix $X\in
\cSE _B$ are nilpotent (see the matrix $A$ of Example 1).
\begin{lemma}\label{1'} For $\, i,i'\in \{ 1,\ldots ,u\} \, $,
$\, j\in \{ q_i-q_{i-1},\ldots ,1\} \, $ and $\ j'\in \{
q_{i'}-q_{i'-1}, \ldots ,1\} \ $
 the maximum rank of  $X_{q_i-j+1,q_{i'}-j'+1}$ for $\; X\in \cSE
_B$ ($\; X\in \cSN _B\; $) is:
  \begin{itemize} \item[a)] $\ \mu
_{q_{i'}}\ \ $ if $\ \ i<i' \ $ ($\ \mu_ {q_i}
>\mu _{q_{i'}}\ $); \item[b)] $\ \mu _{q_i}-1\
\ $ if $\ \ i=i'\ \ $ and $\ \ j\leq j'\ $; \item[c)] $\ \mu
_{q_i}\ \ $ if $\ \ i>i'\ $ or $\; $ if $\ \ i=i'\ $ , $\ j>j'\
$.\end{itemize} \end{lemma} \pf The claim is a consequence of the
fact that if $X\in \cU _B$ then $X\in \cSE _B$ iff $X(i,l)$ is
strictly lower triangular for $i=1,\ldots ,u$ (if $X\in \cSN _B$
then $X(i,l)=X(i,l+1)$ for $l= 1,\ldots , \mu _{q_i}-1 $).
\hspace{4mm} $\square $ \vspace{2mm}
\newline If $X$ is an endomorphism of $K^n$ and $\Lambda $ is a
basis of $K^n$ we will denote by $R _{X, \Lambda }$ the relation
in the set of the elements of $\Lambda $ defined as follows: $\,
w'\ R _{X, \Lambda } \ w\, $ iff $\; X\, w'\, $ has nonzero entry
with respect to $w$. By the form of the matrices of $\cU _B$ and
lemma \ref{1'} we get the following description of the elements of
$\cSE _B$.
\begin{corollary}\label{cor'} There exists a non-empty
open subset of $\, \cSE _B\, $ ($\, \cSN _B\, $) such that if $X$
belongs to it and $\ v_{\mu _{q_i},j}^{l}\ ,\ v_{\mu _{q_{i'}},
j'}^{l'}\in \Delta _B\ $ then $\, v_{\mu _{q_{i'}},j'}^{l'}\ R
_{X, \Delta _B} \ \ v_{\mu _{q_i},j}^l\ $ iff one of the following
conditions holds:
\begin{itemize} \item[$\iota _1 $)] $\, i<i'\, $  and $\,  \mu _{q_i}-l\leq \mu _{q_{i'}}-l'$; \item[$\iota  _2$)]
$ \, i=i'\, $, $\, j\geq j' \, $  and $\, l>l'$ (that is $\, \mu
_{q_i}-l<\mu _{q_{i'}}-l'\, $); \item[$\iota _3 $)] $\, i=i'\ $,
$\, j<j'\, $ and $\, l\geq l'$ (that is $\, \mu _{q_i}-l\leq \mu
_{q_{i'}}-l'\, $);
 \item[$\iota _4$)] $\ i>i'\ $
  and $\ l\geq l'\ $.\end{itemize} \end{corollary} The subalgebra $\cSE _B$ ($\
\cSN _B$) is a maximal nilpotent subalgebra of $\cU _B$ ($\cC
_B$); hence there exists a basis of $K^n$ with respect to which
all the elements of $\cSE _B\ $ ($\ \cSN _B$)
 are  upper triangular. This basis is just $\Delta _B$ with the following new
 order.
We set $ v_{\mu _{q_i},j}^{l}\, \prec \, v_{\mu _{q_{i'}},
j'}^{l'} $ if one of the following conditions holds:
\begin{itemize} \item[$e_1$)] $ \ \mu _{q_i}-l<\mu _{q_{i'}}-l'\ $; \item[$e_2$)] $\ \mu _{q_i}-l=\mu _{q_{i'}}-l'\
$ and $\ i<i'\ $ (hence $l>l'$); \item[$e_3$)] $\ \mu _{q_i}-l=\mu
_{q_{i'}}-l'\ $, $\ i=i'\ $ (hence $l=l'$) and $\ j>j'\ $.
\end{itemize} Let $\Delta _{B,\prec }$ be the basis of $K^n$ which
has the same elements as $\Delta _B$ with the order $\prec $. Then
by Corollary \ref{cor'} we have that the representation of all the
elements of $\cSE _B\ $ ($\ \cSN _B$) with respect to $\Delta
_{B,\prec }$ is upper triangular. \vspace{2mm}
\newline
{\bf Example 2} $\ $ Let $n=13$ and $B=(4,3^2,2,1)$. In this case
we have:
$$\Delta _B=\{ \underbrace{v_{4,1}^4, v_{4,1}^3,
v_{4,1}^2, v_{4,1}^1}, \underbrace{v_{3,2}^3, v_{3,2}^2,
v_{3,2}^1},\underbrace{ v_{3,1}^3, v_{3,1}^2,
v_{3,1}^1},\underbrace{ v_{2,1}^2, v_{2,1}^1},\underbrace{
v_{1,1}^1}\} \ , $$ $$ \Delta _{B,\prec }=\{
\underbrace{v_{4,1}^4, v_{3,1}^3, v_{3,2}^3, v_{2,1}^2,
v_{1,1}^1},\underbrace{ v_{4,1}^3, v_{3,1}^2, v_{3,2}^2,
v_{2,1}^1},\underbrace{ v_{4,1}^2, v_{3,1}^1, v_{3,2}^1},
\underbrace{v_{4,1}^1}\} \ .$$ The matrices of an
 endomorphism
$\overline A$ of $K^{13}$ with respect to $\Delta _B$ and with
respect to $\Delta _{B,\prec } $ are respectively as
follows:\vspace{2mm} {\small{$$ \pmatrix{ 0 & a & b & c & | & p &
q & r & | & v & z & w & | & \alpha & \beta & | & \lambda \cr & 0 &
a & b & | & & p & q & | & & v & z & | &   & \alpha & | & 0 \cr & &
0 & a & | & & & p & | & & & v & | & &
 & | & \cr  & & & 0 & | & & & & | & & & & | & &  & | & &
\cr - & -& -& -& -& -& -& -& -& -& - & -& -& -& -& - \cr & s & t &
u & | & 0 & i & l & | & 0 & m & n & | & \rho & \sigma & | & \pi
\cr & & s & t & | & & 0 & i & | & & 0 &
 m & | &  & \rho & | & 0 \cr & & & s & | & & & 0 & | & & & 0
& | & & & \cr - & -& -& -& -& -& -& -& -& -& - & -& -& -& -& - \cr
& j & y & k & | & d & e & f & | & 0 & g & h & | & \xi & \zeta & |
& \theta \cr & & j & y & | & & d & e & | &  & 0 & g & | & & \xi &
| & 0 \cr & & & j & | & & & d & | & & & 0 & | & & & | & \cr - & -&
-& -& -& -& -& -& -& -& - & -& -& -& -& - \cr & {\bf 0} & \delta &
\epsilon & | & & \eta & \nu & | & & \tau & \gamma & | & 0 & o & |
& \phi \cr & & {\bf 0} & \delta & | & & & \eta & | & & &  \tau & |
& & 0 & | & \cr - & -& -& -& -& -& -& -& -& -& - & -& -& -& -& -
\cr & {\bf 0} & {\bf 0} & \omega & | & & {\bf 0} & \psi & | & &
{\bf 0} & \iota & | & & \chi & | & 0 \cr }\ .$$}}\vspace{2mm}
{\small {$$\pmatrix{0 & v & p & \alpha & \lambda & | & a & z & q &
\beta & | & b & w & r & | & c \cr & 0 & d & \xi & \theta & | & j &
g & e & \zeta & | & y & h & f & | & k \cr & & 0 & \rho & \pi &  |
& s & m & i & \sigma & | & t & n & l & | & u \cr & & & 0 & \phi &
| & {\bf 0} & \tau & \eta & o & | & \delta  & \gamma  & \nu & | &
\epsilon \cr  & & & & 0 & | & {\bf 0} & {\bf 0} & {\bf 0} & \chi &
| & {\bf 0} & \iota & \psi & | & \omega \cr - & -& -& -& -& -& -&
-& -& -& - & -& -& -& -& - \cr & & & & & | & 0 & v & p & \alpha &
| & a & z & q & | & b \cr & & & & & | & & 0 & d & \xi & | & j & g
& e & | & y \cr &  & & & & | & & & 0 & \rho & | & s & m & i & | &
 t \cr &  & & & & |&  & & & 0 & | & {\bf 0} & \tau & \eta &  | &
 \delta \cr - & -& -& -& -& -& -& -& -& -& - & -& -& -& -& - \cr &
 & & & &  | & & & & & | & 0 & v & p & | & a \cr &
& & & &  | & & & &  & | &  & 0 & d & | & j \cr &  & & & & | & & &
& & | & & & 0 & | & s \cr  - & -& -& -& -& -& -& -& -& -& - & -&
-& -& -& - \cr &  & & & & | & & & & & | & & & & | & 0 \cr  } \
.$$}}  In this section and in the following one we will represent
any endomorphism of $K^n$ with respect to the basis $\Delta
_{B,\prec }\ $.  For $\ h=0,\ldots , \mu _{q_1}-1\ $ let
$$\Delta _{B, h}=\{ v_{\mu _{q_i},j}^l\ \in \ \Delta _B\ | \ j=q_i -q_i-1, \ldots ,1, \ \mu
_{q_i}-l=h\, \} \ ,$$ with the order induced by $\prec \ $. We
have that $|\Delta _{B, h}|\, \geq \, |\Delta _{B, h'}|\ $ if $\,
h<h'$. Let $\pi _h$ be the canonical projection of $K^n$ onto
$\langle \Delta _{B, h}\rangle $; for $X\in M(n,K)\ $ and $\
h,k\in \{ 0, \ldots , \mu _{q_1}-1\} \ $ let $ \ X_{h,k}=\pi _h\,
\circ \, X|_{\langle \Delta _{B, k} \rangle }\ $; we consider $X$
as a block matrix:
$$X=(X_{h,k})\ , \ h,k\in \{ 0, \ldots , \mu _{q_1}-1\} \ .$$ The
following results of this section have the purpose to describe the
representation of the algebra $\cSN _B$ with respect to $\Delta
_{B,\prec }\ $, as it appears in Example 2. \newline
 Let $\, A\in \cSN _B\, $. If we
cancel the column of the entries of $\, v_{\mu _{q_1},q_1}^{l}\, $
and the row of the entries with respect to $\, v_{\mu _{q_1},
q_1}^{l}\, $
 for $\, l=1,\ldots ,\mu _{q_1}\,  $ (that is the first
row of each row of blocks and the first column of each column of
blocks) we get a matrix $\, A_{(1)}\, $ of $\, \cSN _{B_{(1)}}\,
$, where $\, B_{(1)}=(\mu _2,\ldots ,\mu _t)\, $. If we instead
cancel the column of the entries of $\, v_{\mu _{q_i}, j}^1$ and
the row of the entries with respect to $\, v_{\mu _{q_i}, j}^1$
for $\, i=1,\ldots ,u\, $ and $\, j=1, \ldots ,q_{i}-q_{i-1}\, $
we get a matrix $\,  A^{(1)}\, $ of $\, \cSN _{B^{(1)}}\, $, where
$\, B^{(1)}=(\mu _{1}-1,\ldots ,\mu _t-1)\, $ (here we cancel the
$0$'s).\newline We observe that for $X\in M(n,K)$ we have $X\in
\cSN _B$ if and only if $X$ is a strictly upper triangular matrix
with the following property: any block $X_{h,k}$ for $h,k=0,\ldots
, \mu _{ q_1}-2$ has the block $X_{h+1,k+1}$ as a submatrix and if
they are different then $X_{h,k}$ is obtained by writing $0$'s
under $X_{h+1,k+1}$ and (or) other columns on its right. This is
expressed by the following Proposition.
\begin{proposition}\label{forma} With respect to $\Delta _{B,\prec }$ the
variety $\cSN _B$ is the affine space of all the strictly upper
triangular matrices $A$ such that:\begin{itemize} \item[i)] if $\,
\mu _{q_{i'}}-\mu _{q_i}>l'-l>0\, $ (hence $\mu _{q_{i'}}-\mu
_{q_i}>1$) the entry of $\, A\, v_{\mu _{q_{i'}},j'}^{l'}\, $ with
respect to $\, v_{\mu _{q_i},j}^l\, $ is $0$; \item[ii)] for $\,
k\in \{ 1,\ldots ,\mu _{q_1}-1\} \, $ and $\, h\in \{ 1,\ldots
,k\} \, $ the entry of $\, A_{h,k}\, $ of indices $\, (i,j)\, $ is
equal to the entry of $\, A_{h-1,k-1}\, $ of indices $\, (i,j)\,
$.
\end{itemize} \end{proposition}
\pf The claim i) can be proved by Corollary \ref{cor'}; the claim
ii) can be proved by lemma \ref{2}, Corollary \ref{cor'} and
induction on $n$ (we can consider $A_{(1)}$ or $ A^{(1)}$).
\hspace{4mm} $\square $\vspace{2mm}
\newline
We will denote by $X=(x_{i,j})$ the generic element of $\cSE _B$,
whose nonzero entries are the coordinates of the affine space
$\cSE _B$. For $h=0,\ldots ,\mu _{q_1}-1$ let $B^{(h)}=(\mu
_1-h,\ldots ,\mu _t-h)$
 (omitting the values which are not positive) and let
 $$t_h=\sum _{l=0}^h\Big| \Delta
_{B, l}\Big|\ .$$ Then $t_0=t$; we set $t_{-1}=0$. By Proposition
\ref{forma} we get the following result. \begin{lemma} \label{ind}
If $h\in \{ \mu _{q_1}-1,\ldots ,0\} $ the submatrix of $X$
obtained by choosing the last $n-t_{h-1}$ rows and columns is the
generic element of $\cSE _{B^{(h)}}$.
\end{lemma}
\begin{corollary}\label{cor''} For  $\ h,k\in \{ 0,\ldots ,\mu _{q_1}-1\}
\ $ and  $\ m\in \N \ $  we have that: \begin{itemize} \item[iii)]
if $\ Y\in M(|\Delta _{B, h}|,K)\ $ is strictly upper triangular
there exists $\ A\in \cSN _B\ $ such that $\ (A^m)_{h,h}=Y^m$;
\item[iv)] if the entry of $\, (A^m)_{h,k}\, $ of indices $(i,j)$
is $0$ for all $\, A\in \cSN _B\, $ then, for all $\, A\in \cSN
_B\, $, we have that: \begin{itemize} \item[iv')] the entry of $\,
(A^m)_{h,k}\, $ of indices $ (i',j') $ is $0$ for $\, i'=i,\ldots
,|\Delta _{B,h}|\, $ and $\, j'=1,\ldots ,j\, $; \item[iv'')] the
entry of indices $ (i,j) $ of $\, (A^m)_{h,k-1}\, $ (if $\, k\neq
0\, $) and of $\, (A^m)_{h+1,k} \, $ (if $\, h\neq \mu _{q_1}-1\,
$ and $\, i\leq | \Delta _{B,h+1}|\,$) is also $0$. \end{itemize}
\end{itemize}\end{corollary} \pf We can prove the claim by
Proposition \ref{forma}, using for simplicity induction on $n$ and
lemma \ref{ind}. By ii) of Proposition \ref{forma} in order to
prove iv'') it is enough to prove the following claim: if $\ k\in
\{ 1,\ldots ,\mu _{q_1}-1\} \  $ and the entry of indices $\
(i,j)\ $ of $\ (A^m)_{0,k}\ $ is $\ 0\ $ for all $\ A\in \cSN _B\
$ then the entry of indices $\ (i,j)\ $ of $\ (A^m)_{0,k-1}\ $ is
also $\ 0\ $ for all $\ A\in \cSN _B$. It can be proved by
induction on $\ |\Delta _{B, 0}|-i\ $ (if $\ i=|\Delta _{B, 0}|\ $
the claim is true, since $\ (A_{(1)})^{m-1}\ $ has the property
iv) by the inductive hypothesis). \hspace{4mm} $\square $
\vspace{2mm} \newline We will denote by $A=(a_{i,j})$ the generic
element of $\cSN _B$, whose nonzero different entries are the
coordinates of the affine space $\cSN _B$. By Proposition
\ref{forma} we also get the following result. \begin{corollary}
\label{added} For $i=1,\ldots ,n $ the entry of $A$ of indices
$(i,i+1)$ is different from $0$ iff $\, i\neq t_h$ for $h=0,
\ldots , \mu _{q_1}-1$ or $i=t_h$ and $|\Delta _{B,h}|=|\Delta
_{B,\mu _{q_1}-2}|$. \end{corollary} By Corollary \ref{cor''} we
get the following result.
\begin{corollary} \label{added'} If $i\in \{ 1,\ldots ,n\} $,
$J\subseteq \{ i+2,\ldots ,n\} $ and there exists $l\in J $ such
that $a_{i,l}\neq 0$ then there exists $j\in J $ such that $a_{i,
j}$ is trascendent over the field of the rational functions in the
entries of the set $$\{ a_{r,s}\ | \  r=i, \ldots ,n, \ s\in J,\
(r,s)\neq (i, j)\} \ .$$
\end{corollary}
\pf Let $h_l\in \{ \mu _{q_1} -1, \ldots ,0\} $ be such that
$t_{h_l-1}<l\leq t_{h_l}$;  let $h_J$ be the maximum of the set $$
\Big\{ 0\Big\} \, \cup \, \Big \{ h\in \{ \mu _{q_1} -1, \ldots
,1\} \ | \ J\not \subseteq \{ 1, \ldots , t_{h-1}\, \} \Big\} \
.$$ If $h_J=h_l\, $ by ii) of Proposition \ref{forma} we can set
$j=l$, since $a_{i,l}$ is different from all the other entries of
the submatrix $(a_{r,s})$, $i=1,\ldots ,n$, $s\in J$. Hence we can
prove the claim by induction on $h_J-h_l$. If $a_{i,l}$ is not
different from all the entries of the set $$\{ a_{r,s}\ | \  r=i,
\ldots ,n, \ s\in J,\ (r,s)\neq (i, l)\} $$ then by ii) of
Proposition \ref{forma} there exist $l'\in J\, $ and $\, h_{l'}\in
\{ \mu _{q_1}-1,\ldots ,0\} \, $ such that $\, h_{l'}>h_l\, $, $\;
t_{h_{l'}-1}<l'\leq t_{h_{l'}}\, $, $\;
l'-t_{h_{l'}-1}=l-t_{h_{l}-1}$. But then by iv'') of Corollary
\ref{cor''} we get that $a_{i,l'}\neq 0$. Hence the claim follows
by the inductive hypothesis. \
 \hspace{4mm} $\square $
\section{Properties of $\cSE _B$ and
$\cSN _B$} For any $n\times n$ strictly upper triangular matrix
$Y=(y_{i,j})$ over any ring let $\Phi _Y$ be the map from $\{
1,\ldots ,n\} $ to $\{ 2,\ldots ,n+1\} $ which is defined as
follows:
$$\Phi _Y(i)=\left\{ \begin{array}{ll} n+1 & \mbox{\rm if
$y_{i,l}=0$ for $l=i+1,\ldots ,n$} \\ \\ \mbox{\rm min } \{ l\in
\{ i+1,\ldots ,n\} \, | \, y_{i,l} \neq 0\} & \mbox{\rm otherwise}
\ .\end{array} \right. $$ For any map $\Phi $ from $\{ 1,\ldots
,n\} $ to $\{ 2,\ldots ,n+1\} $ let $N_{\Phi }$ be the set of all
the $n\times n $ matrices $Y$ over $K$ which are strictly upper
triangular and such that $\Phi _Y=\Phi $ (which is a quasi
projective subvariety).
\begin{lemma} \label{phi} If $\; \Phi $ is not decreasing and the
restriction of $\; \Phi $ to $\Phi ^{-1}\big( \{ 2,\ldots ,n\}
\big) $ is increasing there exists a nilpotent orbit which
contains $N_{\Phi }$.
\end{lemma} \pf For $k\in \N $ we can define the map $\Phi ^k$
from $\{ 1,\ldots ,n\} $ to $\{ 2,\ldots ,n+1\} $ by induction on
$k$, as follows: $\Phi ^1=\Phi $ and $$ \Phi ^k(i)=\left\{
\begin{array}{ll} n+1 & \mbox{\rm if $\Phi ^{k-1}(i)=n+1$} \\ \\
\Phi \big( \Phi ^{k-1}(i)\big) & \mbox{\rm otherwise} \
.\end{array} \right. $$ If $\Phi $ is not decreasing and the
restriction of $\Phi $ to $\Phi ^{-1}\big( \{ 2,\ldots ,n\} \big)
$ is increasing then $\Phi ^k$ is not decreasing and the
restriction of $\Phi ^k$ to $\big( \Phi ^k\big) ^{-1}\big( \{
2,\ldots ,n\} \big) $ is increasing for all $k\in \N $; moreover
the rank of all the elements of $N_{\Phi ^k}$ is the number of
their nonzero rows, that is $$\Big| \{ i\in \{1,\ldots ,n\} \ |\
\Phi ^k(i)\neq n+1\} \Big| \ .$$ If $Y\in N_{\Phi } $ then $Y^k\in
N_{\Phi ^k}$ for all $k\in \N $; this can be proved by induction
on $k$ and considering the transpose matrices of $Y$ and
$Y^{k-1}$. Hence we get that rank $Y^k$ is the same for all $Y\in
N_{\Phi }$, which proves the claim.\hspace{4mm}$\square $
\vspace{2mm}
\newline
Let {\tt id} be the identity map in $\{ 1,\ldots ,n\} $ and let
$I_n$ be the identity matrix of order $n$ over $K$. For any
permutation $\sigma $ of $\{1,\ldots ,n\} $ let
 $G_{\sigma }$ be the matrix of the endomorphism of $K^n$ defined
 by $e_i\longmapsto e_{\sigma (i)}$. If $i,j\in \{1,\ldots ,n\} $
 and $i<j$ let $\sigma _{i,j}$ be the permutation
 $(i,j,j-1,j-2,\ldots ,i+1)$; we set $\sigma _{i,i}= {\tt id}$. \newline
If $Y$ is a strictly upper triangular matrix of order $n$ over a
polynomial ring there exists a matrix $Y'$ over the corresponding
field of the rational functions such that:
\begin{itemize} \item[1)] $Y$ is conjugated to $Y'$ by a product of
matrices each one of which is a product of the form $G_{\sigma
_{i,j}} \cdot G$, where $G$ is upper triangular and the entries of
$G$ are rational functions in the entries of $Y$; \item[2)] $\Phi
_{Y'}$ is not decreasing and the restriction of $\Phi _{Y'}$ to
$(\Phi _{Y'}) ^{-1}\big( \{ 2,\ldots ,n\} \big) $ is increasing.
\end{itemize} In the next Example we show how the matrices
$G_{\sigma _{i,j}}$ and $G$ can be found if $Y$ is the generic
element $A$ of $\cSN_B$, in order to prove that we can get a map
$\Phi _{A'}$ which is not dependent from the Toeplitz conditions.
We will denote by $A_{\{ i\} }$ and $A^{ \{ i\} } $ respectively
the row and the column of $A$ of index $i$ for $i=1,\ldots ,n$.
\vspace{2mm}\newline {\bf Example 3} Let $A$ be the second matrix
of Example 2, where the basis is $\{ e_1,e_2,\ldots ,e_{13}\} $.
If we set $e'_{10}=e_{10}+\tau v^{-1}e_9$ the representation of
$A$ with respect to the basis $\{ e_1, e_2,\ldots , e_9, e'_{10},
e_{11}, e_{12}, e_{13}\} $ is obtained from $A$ by replacing
$A_{\{ 9\} }$ with $(0,\ldots , 0, \eta -\tau v^{-1}, \delta -\tau
v^{-1} )$ and $A^{\{ 10\} }$ with $A^{\{ 10\} } +\tau v^{-1} A^{\{
9\} }$. We also set $e'_{11}=e_{11}+(\eta -\tau v^{-1})d^{-1}\,
e_9\, $ and $\, e'_{12}=e_{12}+(\delta -\tau v^{-1}-(\eta -\tau
v^{-1})d^{-1})s^{-1}\, e_9$; let $G_9$ be the matrix of the map
defined by $e_i\longrightarrow e'_i$ for $i=10, 11, 12$ and
$e_i\longrightarrow e_i$ for $i=1,2,\ldots ,9,13$, that is {\small
$\pmatrix{ I_8 & 0 \cr 0 & G \cr }$} where {\small
$$\begin{array}{lr} &
G=\pmatrix{ 1 &  \tau v^{-1} & (\eta -\tau v^{-1})d^{-1} & (\delta
-\tau v^{-1}-(\eta -\tau v^{-1} )d^{-1})s^{-1} & 0 \cr & 1 & 0 & 0
& 0 \cr & & 1 & 0 & 0 \cr & & & 1 & 0 \cr & & & & 1 \cr } \
.\end{array}
$$} The matrix $(G_9)^{-1} AG_9$ is obtained from $A$ by replacing $A_{\{ 9\}
}$ with the zero row and
$$A^{\{ 10\} } \ \mbox{\rm with }\   A^{\{ 10\} }+\tau v^{-1}
A^{\{ 9\} } \ , \quad A^{\{ 11\} } \ \mbox{\rm with }\     A^{\{
11\} }+(\eta -\tau v^{-1})d^{-1}A^{\{ 9\} }\ ,$$
$$A^{\{ 10\} } \ \mbox{\rm with }\    A^{\{ 10\} }+(\delta -\tau v^{-1}-(\eta -\tau
v^{-1})d^{-1})s^{-1} A^{\{ 9\} }\ .$$ If $y$ is an entry of $A$,
$y"$ will denote a rational function in entries of $A$ which
appear in the same row as $y$ but in columns of smaller index;
$y'$ will denote a sum $y+y"$. We have $\sigma
_{9,12}=(9,12,11,10)$ and the matrix $G_{\sigma _{9,12}
}^{-1}(G_9)^{-1}AG_9G_{\sigma _{9,12}}$ is {$$\pmatrix{0 & v & p &
\alpha & \lambda  & a & z & q & b' & w' & r'  & \beta & c \cr & 0
& d & \xi & \theta  & j & g & e   & y' & h' & f'  & \zeta & k \cr
& & 0 & \rho & \pi  & s & m & i   & t' & n' & l'  & \sigma & u \cr
& & & 0 & \phi  & {\bf 0} & \tau & \eta & \delta ' & \gamma ' &
\nu ' & o & \epsilon \cr & & & & 0  & {\bf 0} & {\bf 0} & {\bf 0}
&  \chi " & \iota '& \psi ' & \chi & \omega \cr  & & & &  & 0 & v
& p   & a' & z' & q'  & \alpha & b \cr & & & &  & & 0 & d  & j' &
g' & e'  & \xi & y \cr & & & & & & & 0  & s' & m' & i'
 & \rho &
 t \cr  &
 & & &  & & &  & 0 & v & p  &  0 & a \cr &
& & &   & & &  &  & 0 & d  & 0 & j \cr &  & & &  & & &  & & & 0  &
{\bf 0} & s \cr  & & & & &  & &  &  &  &   & 0 &
 0 \cr  &  & & &  & & &  & & &  &  & 0\cr  } \
.$$} Now we set
$$e'_8=e_8+ \chi "(s')^{-1}e_5\ , \quad e'_{9}=e_{9}+(\iota '-\chi " (s')^{-1})v^{-1}\, e_5\ ,$$
$$e'_{10}=e_{10}+(\psi ' -\chi " (s')^{-1}-(\iota ' -\chi " (s')^{-1})v^{-1})d^{-1}\,
e_5$$ and $G_5$ to be the matrix of the map defined by
$e_i\longrightarrow e'_i$ for $i=8, 9, 10$ and $e_i\longrightarrow
e_i$ for $i=1,2,\ldots ,7,11,12,13$; then the matrix $$\,
(G_5)^{-1}G_{\sigma _{9,12} }^{-1}(G_9)^{-1}AG_9G_{\sigma
_{9,12}}G_5\, $$ has the entries of indices $(5,l)$, $\; l=1,
\ldots , 11$ equal to $0$ while the entry of indices $(5,12)$ is
$$\chi  -\chi " (s')^{-1}-(\iota ' -\chi " (s')^{-1})v^{-1} -(\psi '
-\chi " (s')^{-1}-(\iota ' -\chi "(s')^{-1})v^{-1} )d^{-1}=$$ $$
\chi -\chi "(s')^{-1} (1-v^{-1}-d^{-1}+v^{-1}d^{-1})-\iota '
v^{-1}(1-d^{-1})-\psi ' d^{-1} $$ and will be denoted by $\chi
^{\ast }$.  By corollary \ref{added'} we get that $\chi ^{\ast }$
is not $0$ (in fact the coefficients of the elements $\chi $,
$\chi "$, $\iota '$, $\psi '$ cannot be equal to the zero
polynomial; for example $\psi '$ is a trascendent element as
expressed in that corollary). The entry of indices $(5,13)$ is
obtained by a similar formula and will be denoted by $\omega
^{\ast }$. Let $\sigma _{5,10}=(5,10,9,8,7,6)$. The matrix
$$(G_{\sigma _{5,10}})^{-1}(G_5)^{-1}G_{\sigma _{9,12}
}^{-1}(G_9)^{-1}AG_9G_{\sigma _{9,12}}G_5G_{\sigma _{5,10}}$$ is
{$$\pmatrix{0 & v & p & \alpha  & a & z & q  & b' & w' & \lambda &
r'  & \beta & c \cr & 0 & d & \xi
  & j & g & e   & y' & h' & \theta & f'  & \zeta & k \cr &
& 0 & \rho  & s & m & i   & t' & n'& \pi  & l'  & \sigma & u \cr &
& & 0   & {\bf 0} & \tau & \eta  & \delta ' & \gamma ' & \phi &
\nu ' & o & \epsilon \cr  & & &  & 0 & v & p   & a' & z'& 0 & q' &
\alpha & b \cr & & &  & & 0 & d  & j' & g' & 0& e'  & \xi & y \cr
& & &  & & & 0  & s' & m' & 0 & i'  & \rho &
 t \cr &
 & &  & & &  & 0 & v & 0 & p  &  0 & a \cr &
& &  & & &  &  & 0 & {\bf 0} & d  & 0 & j \cr & & &  &  &  & &  &
& 0 & 0 & \chi ^{\ast }& \omega ^{\ast } \cr & & &  & & &  & & &
& 0 &
 0 & s \cr & & & & & & &  &  &  &   & 0 &
 0 \cr  &  & &   & & &  & & &  &  & & 0\cr  } \
.$$} The next step will be done by the following definitions:
$e'_5= e_5+\tau v^{-1}e_4\, ,$
$$e'_6=e_6+(\eta -\tau v^{-1})d^{-1}e_4\, , \quad e'_7=e_7+(\delta
'-\tau v^{-1}-(\eta -\tau v^{-1})d^{-1})(s')^{-1} e_4\, ,$$
$$e'_8=e_8+(\gamma '-\tau v^{-1}-(\eta -\tau v^{-1})d^{-1}-(\delta
' -\tau v^{-1}-(\eta -\tau v^{-1})d^{-1})(s')^{-1})v^{-1}e_4$$ and
$\sigma _{4,8}=(4,8,7,6,5)$; we will get an entry $\phi ^{\ast }$
of indices $(8,10)$. Going on in this way we get a matrix $A'$
conjugated to $A$ such that $\Phi _{A'} $ is not decreasing and
the restriction of $\Phi _{A'} $ to $(\Phi _{A'} )^{-1} (
\{2,\ldots ,13\} )$ is increasing. If we had considered the
generic element $X$ of $\cSE _B$ instead of $A$, since $X$ has the
same zero entries as $A$ we would have obtained a matrix $X'$
conjugated to $X$ and such that $\Phi _{A'} =\Phi _{X'} $.
\vspace{2mm}
\newline
In the next Proposition we will extend the results shown in
Example 3 to any partition $B$, showing that there exists an open
subset of $\cSE _B$ whose intersection with $\cSN _B$ is not empty
such that all the elements of it are conjugated to elements of
$N_{\Phi }$ for a certain $\Phi $ which is not decreasing and such
that the restriction of $\Phi $ to $\Phi ^{-1} (\{ 2,\ldots ,n\}
)$ is increasing.
 We set $\{ e_1,\ldots ,e_n\} =\Delta _{B,\prec }\; $(the basis with
 respect to which we represent the elements of $\cSE
 _B$).
We denote by $G_n, G_{n-1}, \ldots ,G_1$ a sequence of $n\times n$
upper triangular matrices over the field of the rational functions
over
 $\cSE _B$, where $G_l$ is such that the elements along the main diagonal are $1$
while the elements which don't belong neither to the main diagonal
nor to the row of index $l$ are $0$, for $l=n,\ldots ,1$. For
 $l=n,\ldots ,1$ and $i=l+1,\ldots ,n$ we set $g _{l,i}$ to be
 the entry of $G_l(e_i)$ with respect to $e_l$; we set $g
 _{l,l}=0$. We will denote by $\sigma _n, \sigma _{n-1}, \ldots ,\sigma _1$
a sequence of permutations of $\{ 1,\ldots ,n\} $. We denote by
$X$ the generic element of $\cSE _B$; we set
$$X^{[n]}=X\, , \qquad X^{[l]}=\big( G_{\sigma _l}\big) ^{-1}
 \big( G_l\big) ^{-1}X^{[l+1]}G_l
G_{\sigma _l}\qquad \mbox{\rm for } \ l=n-1,\ldots ,1$$
 and $X^{[l]}=\big( x_{i,j}^{[l]}\big) $, $i,j=1,\dots ,n$. For $l=n,\ldots ,1$ let $j_l$ be the maximum of the set
 $$\{ k\ \in \ \{ l,\ldots ,n\} \ | \ \Phi _{X^{[l+1]}}(k)\leq
 \Phi _{ (G_l)^{-1}X^{ [l+1]} G_l} (l)\} \ .$$
 For $l=n,\ldots ,1$ the matrix $G_l$ and the permutation $\sigma
 _l$ are inductively defined according to the
 properties $P_{G_l}$ and $P_{\sigma_l}$, which we are going to
 explain.
 \begin{itemize} \item[$P_{G_l}\ )$] If $\, \Phi
 _{X^{[l+1]}}(l)=n+1\, $ or $\, \Phi _{X^{[l+1]}}(l)<\Phi
 _{X^{[l+1]}}(l+1)\, $ then $\, G_l\, $ is $\, I_n\, $, otherwise $g_{l,i}$ is defined
 by induction on $i$ for $i=l+1,\ldots ,n$ in the following way:
if $
 \Phi _{X^{[l+1]}}(i)=\Phi _{X^{[l+1]}}(i-1)+1 $  or $ x_{l,k}^{[l+1]}=
 0$ for all $k\in \{ 2, \ldots , \Phi
 _{X^{[l+1]}}(i)-1\} $ then $$
g_{l,i}= \Big( x_{l,\Phi _{X^{[l+1]}}(i)}^{[l+1]}- {\displaystyle
\sum _{k=l} ^{i-1}}g _{l,k} \Big) \Big( x_{i,\Phi
_{X^{[l+1]}}(i)}^{[l+1]} \Big) ^{-1}\ ;$$ otherwise (that is, for
greater values of $i$) we have $g_{l,i}=0 \ $;
 \item[$P_{\sigma_l}\ )$]  \
$\sigma _l=\sigma _{l,j_l}\ $. \end{itemize} Let $A^{[l]}$ be the
matrix obtained from $X^{[l]}$ by replacing any variable which has
the smallest indices in one of the equations of $\cSN _B$ (as
subvariety of $\cSE _B$) with the variable of that equation with
greatest indices (equal to it according to that equation).
\begin{proposition}\label{Toeplitz} For $l=n,\ldots ,1$
there exist an upper triangular $n\times n$ matrix $G_l$ with the
property $P_{G_l}$ and a permutation $\sigma _l$ of $\{ 1,\ldots
,n\} $ with the property $P_{\sigma _l}$ such that:
\begin{itemize} \item[i)] the entries of $A^{[l]}$ are rational functions in the coordinates of the affine space $\cSN _B$;
\item[ii)] there exists a map $\Phi _l$ from $\{ 1,\ldots ,n\} $
to $\{ 2,\ldots ,n+1\} $ such that the restriction of $\Phi _l$ to
$\{ l,\ldots ,n\} $ is not decreasing, the restriction of $\Phi
_l$ to $ \Phi _l ^{-1} \Big( \{ l+1,\ldots , n\} \Big)\; \cap \;
\{ l,\ldots ,n\} $ is increasing and $\Phi _{X^{[l]}}=\Phi
_{A^{[l]}}=\Phi _l$.
\end{itemize}
\end{proposition}
\pf  If $\mu _1-\mu _t\leq 1$ then the rank of $X$ and $A$ is
$n-1$, hence we set $G_l= I_n$, $\; \sigma _l=$ {\tt id} for
$l=n,\ldots ,1$. Since the claim is true if $\mu _1-\mu _t\leq 1$
we can prove it by induction on $n$. By the inductive hypothesis
the claim is true for $B^{(1)}$, hence by lemma \ref{ind} we can
assume that it is true for $l\in \{ n,\ldots ,t_0+1\} $ and we
prove it for $l=t_0,\ldots ,1$. Let $G_{t_0}$ be as described in
the property $P_{G_{t_0}}$. By the definition of $G_{t_0}$ we get
that either $\Phi _{(G_{t_0})^{-1}X^{[t_0+1]} G_{t_0}}\,
(t_0)=n+1$ or $\Phi _{(G_{t_0})^{-1}X^{[t_0+1]} G_{t_0}}\, (t_0)$
does not belong to the image of the restriction of $\Phi
_{X^{[t_0+1]}} $ to $\{ t_0+1, \ldots ,n\} $. Still by the
definition of $G_{t_0}$, by Corollary \ref{added'} and by the
inductive hypothesis on $X^{[t_0+1]}$ we get that the matrix which
is obtained from $(G_{t_0})^{-1}X^{[t_0+1]} G_{t_0}$ by replacing
any variable which in one of the equations of $\cSN _B$ has the
smallest indices with the variable of that equation with greatest
indices has the following property: its entries are rational
functions in the coordinates of the affine space $\cSN _B$ and if
$\Phi _{(G_{t_0})^{-1}X^{[t_0+1]} G_{t_0}}\, (t_0)\neq n+1$ the
entry of indices $(t_0,\Phi _{(G_{t_0})^{-1}X^{[t_0+1]} G_{t_0}}\,
(t_0))$ is different from $0$ (this is a consequence of Corollary
\ref{added'}, as it was noticed for the element $\chi ^{\ast }$ in
Example 3). By the property $P_{\sigma_{t_0}}$ and by the
inductive hypothesis on $X^{[t_0+1]}$ we get that $\ X^{[t_0]}$
has the required properties. We can repeat the same argument for
$l=t_0-1,\ldots ,1$. \hspace{4mm} $\square $ \vspace{2mm}
\newline
 By lemma \ref{phi} and Proposition \ref{Toeplitz} we get the following
 result, which claims
the independence of $Q(B)$ from the Toeplitz conditions.
\begin{corollary}\label{t3}  The maximum orbit intersecting
$\cSE _B$ ($\cE _B$) is $\, Q(B)$.
\end{corollary}
\pf Let $\Phi _1$ be as described in Proposition \ref{Toeplitz};
then there exists an open subset $\cA $ of $\cSE _B$ whose
intersection with $\cSN _B$ is not empty such that all the
elements of $\cA $ are conjugated to elements of $N_{\Phi _1}$,
which proves the claim by ii) of Proposition \ref{Toeplitz} and
lemma \ref{phi}. \hspace{4mm} $\square $
\section{The graph associated to $B$}
Let $R_B$ be the relation in the set of the elements of $\Delta
_B$ defined as follows:
$$v_{\mu _{q_{i'}}, j'}^{l'}\ R _B \ v_{\mu
_{q_i},j}^l \ \Longleftrightarrow \   \iota _1) \ \mbox{\rm or }
\iota _2) \ \mbox{\rm or } \iota _3) \ \mbox{\rm or } \iota _4)\
\mbox{\rm of Corollary \ref{cor'} holds} \ .
$$ \begin{proposition} \label{order} The relation $R _B$ in the
set of the elements of $\Delta _B$ is a strict partial order.
\end{proposition} \pf The relation $R _B $ is obviously
antisymmetric; the condition $\iota _1)$ implies $l>l'$, hence it
is also transitive. \hspace{4mm} $\square $\vspace{2mm} \newline
The relation $R_B$ describes a generic element of $\cSE_B$, hence
by Corollary \ref{t3} we get the following result.
\begin{corollary}\label{rel} The maximum nilpotent orbit of the
elements of $\cSN _B$ is determined by the relation $R _B$.
\end{corollary} We will write the vertices of the graph of $R_B$
forming a table as follows. The indices of the rows are the
elements of $\N \cup \{ 0\} $ and the columns have as indices
$\mu _{q_u},\ldots ,\mu _{q_1}$; the element $v_{\mu _{q_i},j}^l$ of $\Delta _B$ is
written in the $\mu _{q_i}-$th column and in the row whose index is
the maximum number $m$ such that there exist elements of $\Delta
_B$ whose images under $X^m$ have nonzero entry with respect to
$v_{\mu _{q_i},j}^l$ for some $X\in \cSE _B$. The
 graph of $R _B$ can be obtained by writing arrows on this table
 according to Corollary \ref{cor'}. We will say that $v_{\mu _{q_{i'}}, j'}^{l'}$
 precedes $v_{\mu
_{q_{i}}, j}^{l}$ in the graph of $R_B$ if the row of $v_{\mu
_{q_{i'}}, j'}^{l'}$ is smaller than the row of $v_{\mu _{q_{i}},
j}^{l}$ (this happens if there is an arrow from $v_{\mu _{q_{i'}},
j'}^{l'}$ to $v_{\mu _{q_{i}}, j}^{l}$). \begin{lemma}
\label{graph} The vectors of $\Delta _B$ are written in the graph
of $R_B$ according to the following rules: \begin{itemize}
\item[a)] $v_{\mu _{q_{i}}, j'}^{l'}$ precedes $v_{\mu _{q_{i}},
j}^{l}$ in the graph of $R_B$ if $l'<l$ or if $l'=l$ and $j'>j$;
\item[b)] if $\mu _{q_{i'}}<\mu _{q_i}$ and $\mu _{q_{i'}}-l'\geq
\mu _{q_i}-l$ then $v_{\mu _{q_{i'}}, j'}^{l'}$ precedes $v_{\mu
_{q_{i}}, j}^{l}$ in the graph of $R_B$; \item[c)] if $\mu
_{q_{i'}}>\mu _{q_i}$ and $l'\leq l$ then $v_{\mu _{q_{i'}},
j'}^{l'}$ precedes $v_{\mu _{q_{i}}, j}^{l}$ in the graph of
$R_B$.
\end{itemize} \end{lemma}\pf The claim a) follows by $\iota _2$)
and $\iota _3$) of Corollary \ref{cor'}; the claims b) and c)
follow respectively from $\iota _1$) and $\iota _3$) of the same
Corollary.\hspace{5mm} $\square $ \vspace{2mm}\newline {\bf
Example 4}\hspace{4mm} For $B=(7,5,2)\, $, $\, B=(2^2,1)\, $, $\,
B=(4,2^2,1)\, $,  $\, B=(2^3)\, $, $\, B=(5,2^3) \,$,  $\,
B=(6,5,2^3)\, $, $\, B=(3^2,2,1)\, $, $\, B=(8^2,6^4, 3^2,2,1)\,
$, $\, B=(4,3^2,2,1)\, $,
 $\,
B=(5,4,3^2,2,1)\, $ and $\, B=(17,15,13,5,4,3^2,2,1)\,$  we
respectively get the following graphs (where we omit the arrows
and we put $\circ   $ in front of the elements of $\Delta ^{\circ
}_B$):
$$ \begin{array}{cccc} & \mbox{\bf 2} & \mbox{\bf 5}& \mbox{\bf 7}
\\ \mbox{\bf 0} & & &  \circ \, v_{71}^1\\ \mbox{\bf 1} & & v_{51}^1 &
\circ \, v_{71}^2  \\ \mbox{\bf 2} & v_{21}^1  & v_{51}^2 & \circ \, v_{71}^3\\
\mbox{\bf 3} & v_{21}^2 & v_{51}^3 & \circ \, v_{71}^4 \\ \mbox{\bf 4} &
& v_{51}^4 & \circ \, v_{71}^5\\ \mbox{\bf 5} &  & v_{51}^5 & \circ \, v_{71}^6\\
\mbox{\bf 6} & & & \circ \, v_{71}^7 \end{array}\qquad \quad
\begin{array}{ccc} & \mbox{\bf 1} & \mbox{\bf 2} \\ \mbox{\bf 0} &
&\circ \, v_{22}^1
\\  \mbox{\bf 1} &  & \circ \, v_{21}^1
\\ \mbox{\bf 2} & \circ \, v_{11}^2&  \\ \mbox{\bf 3} &  & \circ \, v_{22}^2 \\ \mbox{\bf 4} &
&\circ \, v_{21}^2
 \end{array} \qquad \quad  \begin{array}{cccc} & \mbox{\bf 1} & \mbox{\bf 2}
& \mbox{\bf 4}\\ \mbox{\bf 0} &  &  & \circ \, v_{41}^1\\ \mbox{\bf 1} &
 & \circ \, v_{22}^1 & v_{41}^2
\\ \mbox{\bf 2} & & \circ \, v_{21}^1 & \\ \mbox{\bf 3} & \circ \, v_{11}^1& &
v_{41}^3  \\ \mbox{\bf 4} & & \circ \, v_{22}^2 & \\ \mbox{\bf 5} & & \circ \, v_{21}^2 & \\
\mbox{\bf 6} &  & & \circ \, v_{41}^4
\end{array} $$
$$ \begin{array}{c}\begin{array}{cc} & \mbox{\bf 2}\\ \mbox{\bf 0} & \circ \, v_{23}^1 \\ \mbox{\bf
1}& \circ \, v_{22}^1 \\ \mbox{\bf 2} & \circ \, v_{21}^1 \\ \mbox{\bf 3} & \circ \, v_{23}^2
 \\ \mbox{\bf 4} & \circ \, v_{22}^2 \\ \mbox{\bf 5} &
 \circ \, v_{21}^2\end{array}\qquad
\begin{array}{ccc} & \mbox{\bf 2}& \mbox{\bf 5}\\ \mbox{\bf 0} & &\circ \, v_{51}^1  \\ \mbox{\bf
1}&\circ \, v_{23}^1 & v_{51}^2 \\ \mbox{\bf 2} &\circ \, v_{22}^1  &v_{51}^3  \\
\mbox{\bf 3} &\circ \, v_{21}^1 &
 \\ \mbox{\bf 4} &\circ \,  v_{23}^2 &v_{51}^4  \\ \mbox{\bf 5} & \circ \, v_{22}^2 & \\ \mbox{\bf 6} &\circ \, v_{21}^2 &  \\
 \mbox{\bf 7} &  & \circ \, v_{51}^5\end{array} \\ \\ \\ \\ \begin{array}{cccc}  & \mbox{\bf 2} & \mbox{\bf 5} & \mbox{\bf 6} \\
\mbox{\bf 0} & & &\circ \, v_{61}^1   \\
\mbox{\bf 1} & & \circ \, v_{51}^1  &  \\
\mbox{\bf 2} & v_{23}^1 & & \circ \, v_{61}^2 \\
\mbox{\bf 3} & v_{22}^1 & \circ \, v_{51}^2  &  \\
\mbox{\bf 4} & v_{21}^1 & &  \circ \, v_{61}^3  \\
\mbox{\bf 5} & v_{23}^2  & \circ \, v_{51}^3 &   \\
\mbox{\bf 6} & v_{22}^2 &  & \circ \, v_{61}^4
\\ \mbox{\bf 7} &v_{21}^2  & \circ \, v_{51} ^4   &  \\
\mbox{\bf 8} &  &  & \circ \, v_{61}^5 \\ \mbox{\bf 9} &    &\circ \, v_{51}^5  & \\
\mbox{\bf 10} &  & &\circ \,  v_{61}^6
\end{array}\\ \\ \\ \\ \begin{array}{cccc} & \mbox{\bf 1} & \mbox{\bf 2} &
\mbox{\bf 3}\\ \mbox{\bf 0} &  &  & \circ \, v_{32}^1 \\ \mbox{\bf 1} &
 & & \circ \, v_{31}^1
\\ \mbox{\bf 2} & & \circ \, v_{21}^1 & \\ \mbox{\bf 3} & v_{11}^1  & & \circ \, v_{32}^2
\\ \mbox{\bf 4} & & &\circ \,  v_{31}^2 \\ \mbox{\bf 5} & &
\circ \, v_{21}^2 & \\ \mbox{\bf 6} &  & & \circ \, v_{32}^3 \\ \mbox{\bf 7} &
 & & \circ \, v_{31}^3
\end{array}\end{array} \qquad \quad  \begin{array}{cccccc}  & \mbox{\bf 1} & \mbox{\bf 2} & \mbox{\bf 3} & \mbox{\bf 6} & \mbox{\bf 8}\\
\mbox{\bf 0} &  & & & & \circ \, v_{82}^1\\
\mbox{\bf 1} &  & & & & \circ \, v_{81}^1 \\
\mbox{\bf 2} &  & & & \circ \, v_{64}^1 & v_{82}^2 \\
\mbox{\bf 3} &  & & & \circ \, v_{63}^1 & v_{81}^2\\
\mbox{\bf 4} &  & & & \circ \, v_{62}^1 & \\
\mbox{\bf 5} &  & & & \circ \, v_{61}^1 & \\
\mbox{\bf 6} &  & & v_{32}^1  & \circ \, v_{64}^2 & v_{82}^3\\
\mbox{\bf 7} &  &  & v_{31}^1 & \circ \, v_{63}^2 & v_{81}^3 \\
\mbox{\bf 8} &  & v_{21}^1 &  & \circ \, v_{62}^2 &
\\ \mbox{\bf 9} &v_{11}^1  & &  & \circ \,  v_{61} ^2 & \\
\mbox{\bf 10} &  & &v_{32}^2  & \circ \, v_{64}^3 & v_{82}^4\\ \mbox{\bf 11} & &  & v_{31}^2 & \circ \, v_{63}^3 & v_{81}^4\\
\mbox{\bf 12} &  & v_{21}^2 &  &\circ \, v_{62}^3 & \\ \mbox{\bf 13} &  &   & & \circ \, v_{61}^3& \\
\mbox{\bf 14} &
&  &v_{32}^3 &\circ \, v_{64} ^4 & v_{82}^5\\ \mbox{\bf 15} &  &  &v_{31}^3 &\circ \, v_{63}^4& v_{81}^5\\ \mbox{\bf 16} &  & & & \circ \, v_{62}^4 & \\
\mbox{\bf 17} &  & & & \circ \, v_{61}^4 & \\ \mbox{\bf 18} &  &  & &
\circ \, v_{64}^5 & v_{82}^6
\\ \mbox{\bf 19} &  & & &\circ \,  v_{63}^5  & v_{81}^6\\ \mbox{\bf 20} & & & & \circ \, v_{62}^5 &  \\
\mbox{\bf 21} &  & & &\circ \, v_{61}^5 &  \\ \mbox{\bf 22} & & & &
\circ \, v_{64}^6 & v_{82}^7
\\ \mbox{\bf 23} &  & & & \circ \, v_{63}^6 & v_{81}^7  \\ \mbox{\bf 24} &  & & & \circ \, v_{62}^6 & \\
\mbox{\bf 25} & & & & \circ \, v_{61}^6 & \\ \mbox{\bf 26} &  & & & &
\circ \, v_{82}^8\\ \mbox{\bf 27} &  & & & & \circ \, v_{81}^8
\end{array}$$
$$\begin{array}{ccccc} & \mbox{\bf 1}& \mbox{\bf 2} & \mbox{\bf 3} & \mbox{\bf 4}\\
\mbox{\bf 0}&  & & &\circ \,  v_{41}^1 \\ \mbox{\bf 1}& & & \circ \, v_{32}^1  & \\ \mbox{\bf 2}& & & \circ \, v_{31}^1 & \\
\mbox{\bf 3}& &v_{21}^1 & & \circ \, v_{41}^2\\\mbox{\bf 4} & v_{11}^1 & &  \circ \, v_{32}^2 & \\
\mbox{\bf 5}
& & & \circ \, v_{31}^2 & \\ \mbox{\bf 6} & &v_{21}^2& &\circ \, v_{41}^3 \\ \mbox{\bf 7} & & &\circ \, v_{32}^3 & \\
\mbox{\bf 8} & & & \circ \, v_{31}^3& \\ \mbox{\bf 9} & & & &\circ
\, v_{41}^4
\end{array}\qquad \qquad \quad \begin{array}{cccccc}  & \mbox{\bf 1}& \mbox{\bf 2}&\mbox{\bf 3} & \mbox{\bf 4} & \mbox{\bf 5}\\
\mbox{\bf 0} &  & & & &\circ \, v_{51}^1 \\
\mbox{\bf 1} & & & &\circ \, v_{41}^1 & \\ \mbox{\bf 2} & & & \circ \, v_{32}^1 & &  v_{51} ^2\\ \mbox{\bf 3} & & &\circ \, v_{31}^1 & & \\
\mbox{\bf 4} & & v_{21}^1& & \circ \, v_{41}^2& \\ \mbox{\bf 5} & v_{11}^1 & & \circ \, v_{32}^2& &v_{51}^3\\
\mbox{\bf 6} &
& & \circ \, v_{31}^2& & \\ \mbox{\bf 7} & & v_{21}^2& &\circ \, v_{41}^3 & \\ \mbox{\bf 8} &  & &\circ \,  v_{32}^3& & v_{51}^4 \\
\mbox{\bf 9} & & & \circ \, v_{31}^3& & \\ \mbox{\bf 10} & &  & &
\circ \, v_{41}^4&
\\ \mbox{\bf 11} &  & & & &\circ \, v_{51}^5
\vspace{2mm}\end{array}$$
$$\begin{array}{ccccccccc}  & \mbox{\bf 1}& \mbox{\bf 2}&\mbox{\bf 3} & \mbox{\bf 4} & \mbox{\bf 5}& \mbox{\bf 13} & \mbox{\bf 15} & \mbox{\bf 17}\\
\mbox{\bf 0} &  & & & & & & & \circ \, v_{171}^1 \\
\mbox{\bf 1} &  & & & & & & \circ \, v_{151}^1 & v_{171}^2 \\
\mbox{\bf 2} &  & & & & & \circ \, v_{131}^1 & v_{151}^2 & v_{171}^3 \\
\mbox{\bf 3} &  & & & & \circ \, v_{51}^1 & v_{131}^2 & v_{151}^3 & v_{171}^4\\
\mbox{\bf 4} & & & &\circ \, v_{41}^1 & & v_{131}^3 & v_{151}^4 &
v_{171}^5\\ \mbox{\bf 5} & & & \circ \, v_{32}^1 & &  v_{51} ^2
& v_{131}^4 & v_{151}^5 & v_{171}^6\\ \mbox{\bf 6} & & &\circ \, v_{31}^1 & & & v_{131}^5 & v_{151}^6 & v_{171}^7\\
\mbox{\bf 7} & & v_{21}^1& & \circ \, v_{41}^2& & v_{131}^6 & v_{151}^7 &
v_{171}^8\\ \mbox{\bf 8} & v_{11}^1 & & \circ \, v_{32}^2& &v_{51}^3
& v_{131}^{7} & v_{151}^8 & v_{171}^{9}\\
\mbox{\bf 9} &
& & \circ \, v_{31}^2& & & v_{131}^8 & v_{151}^{9} & v_{171}^{10}\\ \mbox{\bf 10} & & v_{21}^2& &\circ \, v_{41}^3 & & v_{131}^{9} & v_{151}^{10} & v_{171}^{11}\\
\mbox{\bf 11} &  & & \circ \, v_{32}^3& & v_{51}^4 & v_{131}^{10} & v_{151}^{11} & v_{171}^{12}\\
\mbox{\bf 12} & & & \circ \, v_{31}^3& & & v_{131}^{11} & v_{151}^{12} &
v_{171}^{13}\\ \mbox{\bf 13} & &  & & \circ \, v_{41}^4 & & v_{131}^{12} &
v_{151}^{13} & v_{171}^{14}
\\ \mbox{\bf 14} &  & & & & \circ \, v_{51}^5 & & v_{151}^{14} &
v_{171}^{15} \\ \mbox{\bf 15} &  & & & & &\circ \, v_{131}^{13} & &
v_{171}^{16} \\ \mbox{\bf 16} &  & & & & & &  \circ \, v_{151}^{15} & \\
\mbox{\bf 17} &  & & & & & & & \circ \, v_{171}^{17}
\vspace{2mm}\end{array}$$
 We recall the definitions of $\, (\tilde
i,\tilde {\epsilon })$, $\, \Delta _B^{\circ, 1}$, $\, \Delta
_B^{\circ, 2}$, $\, \Delta _B^{\circ, 3}$, $\, \Delta _B^{\circ
}$, $\, \widehat B$, $\, Q(B)=(\omega_1, \ldots ,\omega_z)\, $
which were given in section 1; the index of the last row of the
graph of $R_B $ where there are written some vectors is
$\omega_1-1$ and, by Corollary \ref{t3}, the maximum partition
which is associated to elements of $\cSE _B$ is $Q(B)$.\newline By
the definition of the rows of the graph of $R_B$ any row of index
less than $\omega _1$ has at least one vector; by $\iota _2$) and
$\iota _3$) of Corollary \ref{cor'} the first vector and the last
vector of the $\mu _{q_i}-$th  column are respectively $v_{\mu
_{q_i},q_i-q_{i-1}}^1$ and  $v_{\mu _{q_i},1}^{\mu _{q_i}}$ for
$i=1,\ldots ,u$. Moreover we can observe the following less
obvious properties. We consider the injective map $\varphi \ :\
\Delta _B \longrightarrow \{ 0,\ldots ,\omega_1-1\} $ which
associates to any vector the index of its row in the graph of
$R_B$ and the injective map $\phi \, : \, \{ 0,\ldots ,\omega
_1-1\} \longrightarrow \Delta _ B$ such that $\phi (i)$ is the
first vector of the $i-$th row of the graph of $R_B$.
\begin{lemma} \label{next} The basis $\Delta _B$ with the order
$<$ has the following properties:
\begin{itemize} \item[i)] the restriction of the map $\varphi $ to
$\Delta _B^{\circ }$ is a bijection; \item[ii)] the restriction of
$\, \phi $ to $\{ 0,\ldots ,q_{\tilde i +\tilde{\epsilon }} -1\} $
preserves the order and its image is $\Delta _B^{\circ ,1}$; the
restriction of $\, \phi $ to $\{ \omega _1-q_{\tilde i
+\tilde{\epsilon }}+1,\ldots ,\omega _1\} $ reverses the order and
its image is $\Delta _B^{\circ,3}$; \item[iii)] if
 $\, i\in \{ 1,\ldots ,\tilde i+\tilde {\epsilon }\} $, $j\in \{ q_i-q_{i-1},\ldots ,1 \} $  and  $v_{\mu _{q_{i'}},j'}^{l'}$ is
written in the same row of the graph of $R_B$ as $\, v_{\mu
_{q_i},j}^{\mu _{q_i}}$ then $ i'\in \{1,\ldots ,i-1\} $ and
$l'\in \{ \mu _{q_i}+1, \ldots , \mu _{q_{i'}}-1\} $.\end{itemize}
\end{lemma} \pf We can write the graph of $R_B$ starting with writing
the vectors of $\Delta _B^{\circ }$ according to lemma
\ref{graph}. In fact, by Corollary \ref{cor'} and i) of
Proposition \ref{forma} if $\mu _{q_{i'}}-l'<\mu _{q_i}-l$ there
is no arrow from $v_{\mu _{q_{i'}},j'}^{l'}$ to $v_{\mu
_{q_{i}},j}^{l}$ and the same holds if $\mu _{q_{i'}}-\mu
_{q_i}>l'-l>0$; hence by the maximality of $2q_{\tilde i-1}+\mu
_{q_{\tilde i}}(q_{\tilde i}-q_{\tilde i-1})+\tilde {\epsilon }\,
\mu _{q_{\tilde i+1}}(q_{\tilde i+1}-q_{\tilde i})$ if the vectors
of $\Delta _B$ are written in the graph of $R_B$ according to
lemma \ref{graph} then each vector of $\Delta _B-\Delta_B^{\circ
}$ is in the same row as one of the vectors of $\Delta _B^{\circ
}$, which proves i). If $\mu _{q_i}<\mu _{q_{i'}}$ then $v_{\mu
_{q_{i'}},j'}^{1}$ precedes $v_{\mu _{q_i},j}^{1}$ in the graph of
$R_B$ by c) of lemma \ref{graph}, while $v_{\mu _{q_i},j}^{\mu
_{q_i}}$ precedes $v_{\mu _{q_{i'}},j'}^{\mu _{q_{i'}}}$ by b) of
lemma \ref{graph}, hence we get ii). Moreover for $i=1,\ldots
,\tilde i -1$, $j=q_i-q_{i-1},\ldots ,1$ and $l=\mu _{ q_i}-1,
\ldots ,1$ the index of the row of $v_{\mu _{q_i},j}^{l}$ is
smaller than the index of the row of $v_{\mu
_{q_{i+1}},q_{i+1}-q_i}^{\mu _{q_{i+1}}}$ by $\iota _2$), $\iota
_3$) of Corollary \ref{cor'} (see the graph of $R_B$ in the case
$B=(17,15,13, 5,4,3^2,2,1)$), hence we get iii).  \hspace{4mm}
$\square $\vspace{2mm} \newline By i) of lemma \ref{next} we get
Theorem \ref{O}, which was proved by Polona Oblak in \cite{Obl}.
\section{The maximum partition in $\cN _B$ ($\cE _B$)}
In this section we will prove properties of $\cSE _B$ with
arguments which could be used also for $\cSN _B$, hence the proof
of Theorem \ref{eend} is not really dependent from Corollary
\ref{t3}.\newline Let $X$ be the matrix of an endomorphism of
$K^n$ with respect to the basis $\Delta _B$; let $\cSE^{\star }
_B$ be the open subset of $\, \cSE _B\times K^n\, $ of all the
pairs $(X,v)\in N(n,K)\times K^n$ with the following property:
$$\begin{array}{l}   \mbox{\rm if $m\in \N\cup  \{ 0\} $,
$\, m\leq \omega_1-1$ then $ X^m\, v$ has nonzero entry with respect to} \\
\mbox{\rm the vector of  $\Delta _B^{\circ}$ which is written in
the $m-$th
 row of the "graph of $B$"} \ ;  \end{array}
$$ let $K
_B$ be the subset of $K^n$ of all the vectors which have nonzero
entry with respect to $v_{ \mu _{q_1},q_1}^1$.
\begin{proposition} \label{f1}
If $v\in K^n$ the fiber of $v$ with respect to the canonical
projection of $\cSE^{\star } _B$ into $K^n$ is not empty iff $v\in
K _B$.
\end{proposition} \pf By Corollary \ref{cor'} there exists a
non-empty open subset of $\cSE _B$ ($\cSE^{\star } _B$) such that
if $X$ belongs to it the vector $X\, v_{ \mu _{q_1},q_1}^1$ has
nonzero entry with respect to all the elements of $\Delta _B- \{
v_{ \mu _{q_1},q_1}^1\} $. Hence if $v\in K^n$ has nonzero entry
with respect to $v_{ \mu _{q_1},q_1}^1$ and $m,i\in \N \, \cup \,
\{ 0\} $, $ i\geq m$ the condition that $X^m\, v$ has zero entry
with respect to a vector of the $i-$th row of the graph of $R_B$
isn't an identity with respect to $X$. \hspace{4mm} $\square
$\vspace{2mm}
\newline
We recall the following classical lemma (see \cite{Her}).
\begin{lemma} \label{H} If $\, X\in N(n,K)$ has index of nilpotency
$\omega $ and $v\in K^n $ is such that $X^{\omega -1}v\neq 0$
there exists a subspace $V$ of $K^n$ such that $X(V)\subseteq V$
and
$$V \; \oplus \; \langle v, Xv, \ldots , X^{\omega -1} v\rangle \, =
K^n\ .$$ \end{lemma} For $(X,v)\in \cSE  _B\, \times \, K^n$  let
$W_{X,v,0}=\{ 0\} $ and, for $i\in \N $, let $W_{X,v,i}=\langle
X^hJ^kv\, |\, h,k\in \N \cup \{ 0\} \, ,\, k<i \rangle $. The
subspace $W_{X,v,i}$ is stable with respect to $X$;  let $X_{v,i}$ be
the endomorphism of $\; {\displaystyle K^n/ W_{X,v,i}}\; $ defined
by $X_i(W_{X,v,i}+w)=W_{X,v,i}+Xw$.\newline
By Lemma \ref{H} we get the following result.
\begin{proposition} \label{j2} If $(X,v)\in \cSE _B\times K^n$ and $X^{\omega _1-1}v\neq
0$ the partition of $X$ is
 $(\omega_1, \nu _2, \nu _3, \ldots ,\nu _p)$ iff
the partition of $X_{v,1}$ is $(\nu _2, \nu _3, \ldots ,\nu _p)$.
\end{proposition}
 If $(X,v)\in N(n,K)\times K^n$
and $ K^n=\langle X^hJ^kv,\ h,k=1,\ldots ,n\rangle $ one says that
$v$ is cyclic for $X$.
The following Proposition explains a
generalization of a known result for the elements of $\cN _B$.
\begin{lemma}\label{j1} For all $v\in K _B$ the subset of
 $\cSE _B$ of all $X$ such
that $v$ is cyclic for $X$ is not empty.
\end{lemma} \pf In \cite{Na} it has been proved that the
subset of $\cSN _B\times K^n$ of all the pairs $(X,v)$ such that
$v$ is cyclic for $X$ is not empty. The projection of this subset
on $K^n$ is an open subset which is not empty, hence its
intersection with $K _B$ is not empty. Since any element of $K _B$
can be the $\mu_{q_1}-$th element of a Jordan basis for $J$ we get
that $K _B$ is contained in that projection. \hspace{4mm} $\square
$
\begin{proposition}\label{j} If
$(X,v)\in \cSE _B\, \times \, K^n$ and $v$ is cyclic for $X$ then:
\begin{itemize} \item[i)] $X$ has partition $(\nu_1,\ldots ,\nu_{p })$
iff the set $$\Omega _{(X,J)}=\{ \ X^hJ^kv\ |\ ,\ h= \nu _{k+1}-1,\ldots, 0\, , \, k=0,\ldots ,p-1 \ \} $$ is a Jordan basis for
$X$;  \item[ii)] if $X$ has partition $(\nu _1,\ldots ,\nu _{p
})$, $h',k'\in \N \, \cup \, \{ 0\} $ and $k'\geq p $ or $h'\geq
\nu_{k'+1}$ then
$$ X^{h'}J^{k'}v\in \langle \ X^hJ^kv\, |\,
h=\nu _{k+1}-1,\ldots ,0\ ,\ k=0,\ldots ,k'-1\ \rangle \ .$$
\end{itemize} \end{proposition}
\pf  Since $v$ is cyclic for $X$ we have that
$${\displaystyle K^n/W_{X,v,i} }\, =\, \langle \ W_{X,v,i}+X^hJ^kv\ |\ h,k\in
\N \cup \{ 0\} \, ,\ k\geq i\ \rangle \ .$$ Let $X$ have partition
$(\nu _1,\ldots ,\nu _{p })$. Since $\nu _1$ is the index of
nilpotency of $X$ we have $X^{\nu _1-1}v\neq 0$, hence
$W_{X,v,1}=\langle v,Xv,\ldots ,X^{\nu_1-1}v\rangle $. Then $\nu_2$
is the index of nilpotency of $X_{v,1}$ by lemma \ref{H}. Hence
$X^{\nu _2-1}Jv\not \in W_{X,v,1}$. Similarly for $i=2,\ldots ,p -1$
we have that $\ \nu_ {i+1}$ is the index of nilpotency of $X_{v,i}$
and then $X^{\nu _{i+1}-1}J^iv\not \in W_{X,v,i}$.\hspace{4mm}
$\square $
\begin{lemma}\label{directsum} If $\, (X,v)\in \cSE^{\star } _B$ then
$W_{X,v,1}\ \cap \ \langle \Delta _B-\Delta _B^{\circ }\rangle =\{
0\} $. \end{lemma} \pf Let $w$ be an element of $ W_{X,v,1}$
different from $0$;  let $m$ be the minimum element of $ \{
0,\ldots ,\omega_1-1\} $ such that $w$, if represented with respect to
$\{ v, Xv, \ldots , X^{\omega _1 -1} v\} $, has nonzero entry with
respect to $X^mv$.  By the definition of $\cSE ^{\star }_B$ and Corollary \ref{cor'} $w$
 has nonzero entry
with respect to the element of $\Delta _B^{\circ }$ which is
written in the row of index $m$; hence $w\not \in \langle \Delta
_B-\Delta _B^{\circ }\rangle$. \hspace{4mm} $\square $
\vspace{2mm}
\newline If $(X,v)\in \cSE^{\star } _B $ we will denote by $J_{X,v}$ the endomorphism of
${\displaystyle K^n/W_{X,v,1}}$ defined by $
J_{X,v}(W_{X,v,1}+w)=W_{X,v,1}+J(w)$ for $w\in \Delta _B-\Delta
_B^{\circ }$. We observe that if $X\in \cSN _B$ then $J_{X,v}$
commutes with $X_{v,1}$. In general $J_{X,v}$ is not in Jordan canonical
form with respect to the basis $\{ W_{X,v,1}+v_{\mu _{q_i},j}^l\ |\
v_{\mu _{q_i},j}^l\in \Delta _B -\Delta _B^{\circ }\} $; in fact
it is not true that $v_{\mu _{q_i},j}^{ \mu _{q_i}} \in W_{X,v,1} $
for $i = 1,\ldots ,\tilde i-1 $ and $j= q_i-q_{i-1}, \ldots ,1 $
(see the example with $B=(17,15,13, 5,4,3^2,2,1)$, where $\mu
_{q_{\tilde i}} =4$ and $\mu _{q_{\tilde i +\tilde {\epsilon
}}}=3$). In the next Corollary we will show that if $(X,v)\in \cSE^{\star } _B$
the partition of $n-\omega _1$ which corresponds to $J_{X,v}$ is the
partition $\widehat B$ which has been defined in the introduction.
\newline For $\, i=1,\ldots ,\tilde i-1\, $ and $\, j=q_i-q_{i-1},
\ldots ,1\, $ let $\, H(i,j)\, $ be the set of all $\,
(i',j',l')\, $ such that $\, i'\in \{ 1,\ldots i-1\} \, $, $\,
j'\in \{ q_{i'}-q_{i'-1}, \ldots ,1\} \, $ and $l'\in \{ 2, \ldots
, \mu _{q_{i'}}-\mu _{q_i}+1\} $.
\begin{lemma}   \label{nextt} If $(X,v)\in \cSE^{\star }_B$ then for $i=1,\ldots ,\tilde i-1$, $j=q_i-q_{i-1},\ldots ,1$
there exist $a(i,j,i',j',l')\in K$ for $(i',j',l')\in H(i,j)$ such
that:
$$ v_{\mu _{q_i},j}^{\mu _{q_i}}+\sum _{(i',j',l')\in H(i,j)}a(i,j,i',j',l') \, v_{\mu _{q_{i'}},j'}^{l'+\mu _{q_i}-1} \in W_{X,v,1}\ , $$
$$ v_{\mu _{q_i},j}^{\mu _{q_i}-1}+\sum _{(i',j',l')\in H(i,j)}a(i,j,i',j',l') \, v_{\mu _{q_{i'}},j'}^{l'+\mu _{q_i}-2} \not \in W_{X,v,1}\ . $$
\end{lemma} \pf Let $(X,v)\in \cSE _B^{\star }$. By iii) of lemma \ref{graph} for $i=1,\ldots ,\tilde i-1$ and $j=q_i-q_{i-1},\ldots ,1$
there exist $a\in K-\{ 0\} $ and $a(i,j,i',j',l')\in K$ for
$(i',j',l')\in H(i,j)$ such that $$X^{\varphi (v_{\mu
_{q_i},j}^{\mu _{q_i}})}\, v = a\, v_{\mu _{q_i},j}^{\mu
_{q_i}}+\sum _{(i',j',l')\in H(i,j)}a\cdot a(i,j,i',j',l')\,
v_{\mu _{q_{i'}},j'}^{l'+\mu _{q_i}-1}\ ;$$ this proves the first claim. For $m=0,\ldots, \varphi (v_{\mu
_{q_i},j}^{\mu _{q_i}-1})$ the vector
$X^{m}\, v $ has nonzero entry with respect to the vector of $\Delta _B^{\circ }$ which is written in the row of index $m$ of the "graph of $B$", while it has zero entry with respect to the vectors of $\Delta _B^{\circ }$ which are written in the previous rows; hence $$v_{\mu _{q_i},j}^{\mu _{q_i}-1}\ \not \in \ W_{X,v,1}\,+\,  \Big\langle \,v_{\mu _{q_{i'}},j'}^{l'+\mu _{q_i}-2}\ , (i',j',l')\in H(i,j)\, \Big\rangle $$ which proves the second claim.
 \hspace{5mm} $\square $
\vspace{2mm} \newline In the following we will assume that for
$(X,v)\in \cSE _B^{\star }$, $\, i=1,\ldots ,\tilde i-1$, $j=q_i-q_{i-1},\ldots ,1$ and
$(i',j',l')\in H(i,j)$ the elements $a(i,j,i',j',l')\in K$ have
the property expressed in lemma \ref{nextt}.
\begin{corollary} \label{next'} If $(X,v)\in \cSE^{\star } _B$ and we set $$\hat v _{\mu
_{q_i},j}^l=\left\{
\begin{array}{l}v_{\mu _{q_i},j}^{l} + {\displaystyle \sum _{(i',j',l')\in H(i,j)}a(i,j,i',j',l')\, v_{\mu
_{q_{i'}},j'}^{l'+l-1}}\vspace{2mm}
 \\  \quad  \mbox{\rm if } i\in \{ 1,\ldots ,\tilde i-1\} ,\ j\in \{
q_i-q_{i-1},\ldots ,1\} , \ l\in \{ 1,\ldots ,\mu _{q_i}-1\} \\ \\
 v_{\mu _{q_i},j}^{l}  \vspace{1mm}\\  \quad  \mbox{\rm if } i\in \{ \tilde i+\tilde {\epsilon }+1,
\ldots ,u\} , \  j\in \{ q_i-q_{i-1},\ldots ,1\} , \ l\in \{
1,\ldots ,\mu _{q_i}\} \end{array} \right. $$ the representation of
$J_{X,v}$ with respect to the basis
$$\widehat {\Delta }_{B,X,v}\, = \ \{ W_{X,v,1}+\hat v _{\mu _{q_i},j}^l\ |\ v_{\mu _{q_i},j}^l\in
\Delta _B-\Delta _B^{\circ }\} $$ (with the order induced by
$\Delta _B$) has Jordan canonical form and its partition is $\widehat B$.
\end{corollary} \pf  In lemma \ref{nextt} for $i=1,\ldots ,\tilde
i-1$ and $j=q_i-q_{i-1}, \ldots ,1$ we have proved the following
claims:
$$ J^{\mu _{q_i}-1}\, \Big( v_{\mu _{q_i},j}^{1}+\sum
_{(i',j',l')\in H(i,j)}a(i,j,i',j',l') \, v_{\mu
_{q_{i'}},j'}^{l'}\, \Big)  \in W_{X,v,1}\ , $$
$$ J^{\mu _{q_i}-2} \Big( v_{\mu _{q_i},j}^{1}+\sum _{(i',j',l')\in H(i,j)}a(i,j,i',j',l')\, v_{\mu _{q_{i'}},j'}^{l'} \, \Big)
\not \in W_{X,v,1}\ .
$$   By the definition of $J_{X,v}$
they imply that $(J_{X,v})^{\mu _{q_i}-2}\, \Big(W_{X,v,1}+\hat v_{\mu
_{q_i},j}^2\Big)\, =\, W_{X,v,1}$, $\ (J_{X,v})^{\mu _{q_i}-3}\,
\Big(W_{X,v,1}+\hat v_{\mu _{q_i},j}^2\Big) \, \neq  \, W_{X,v,1}$.
Moreover for $i=\tilde i +\tilde {\epsilon }+1, \ldots ,u$ and
$j\in \{ q_i-q_{i-1},\ldots ,1\}$ we have that $\hat v_{\mu
_{q_i},j}^{\mu _{q_i}}=v_{\mu _{q_i},j}^{\mu _{q_i}}\not \in
W_{X,v,1}$, since $v_{\mu _{q_i},j}^{\mu _{q_i}}$ is written in the
same row of the graph of $R_B$ as one of the elements of $\Delta
_B^{\circ }$; this means that $(J_{X,v})^{\mu _{q_i}-1}\, \Big(
W_{X,v,1}+\hat v_{\mu _{q_i},j}^1\Big) \neq W_{X,v,1}$, while
 $(J_{X,v})^{\mu _{q_i}}\, \Big( W_{X,v,1}+\hat v_{\mu
_{q_i},j}^1\Big) = W_{X,v,1}$ since $J^{\mu _{q_i}}\, v_{\mu
_{q_i},j}^1=0$.\hspace{4mm} $\square $ \vspace{2mm}\newline For
any $(X,v)\in \cSE^{\star } _B$ we consider a basis $\widehat
{\Delta }_{B,X,v}$ with the property expressed in lemma \ref{next'}
and we define a map from $\widehat{\Delta } _{B,X,v}$ to $\{
e_1,\ldots ,e_{n-\omega_1} \} $ in the following way: if
$W_{X,v,1}+\hat{v}_{\mu_{q_i},j}^l$ is the $h-$th element of
$\widehat{\Delta }
 _{B,X,v}$ the image of $W_{X,v,1}+\hat{v}_{\mu _{q_i},j}^l$ is $e_h$. This map induces an isomorphism
 $p_{X,v}$ from ${\displaystyle  K^n/
 W_{X,v,1}}\, $ to $\, K^{n-\omega_1}$; we denote by $\pi _B$ the morphism from $\cSE^{\star } _B$ to $N(n-\omega_1,K)$
 defined by $$(X,v)\longmapsto p _{X,v} \,
\circ \, X_{v,1}\, \circ \, \big( p _{X,v}\big) ^{-1}\ .$$ Let $\cQ _B$ be
the open subset of $\cSE^{\star } _B$ of all the pairs $(X,v)$
such that the partition of $X$ is the maximum partition
$Q(B)=(\omega_1,\ldots ,\omega_z)$ which is associated to the
elements of $\cSE _B\, $ ($\cSN _B$).  Let $\widetilde {\cSE }
_{\widehat
 B}$ be the set of all the nilpotent endomorphisms of $K^{n-\omega _1}$ which are conjugated to elements of $\cSE
 _{\widehat B}$.
\begin{proposition} \label{end} The morphism $\pi _B$ has the following properties: \begin{itemize} \item[a)] $\cSE _{\widehat B}\, \subseteq \, \pi _B
\, \big( \cSE^{\star } _B\big)  \ $; \item[b)] $ \pi _B\, \big(
\cQ _B\big)\, \subseteq \, \widetilde {\cSE } _{\widehat
 B} \ $. \end{itemize}
\end{proposition}
\pf \begin{itemize} \item[a)]  We consider the pairs $(X,v)\in
\cSE^{\star }_B$ with the following properties: \begin{itemize}
\item[$P_1)$] there are no arrows from elements of $\Delta
_B^{\circ }$ to elements of $\Delta _B-\Delta _B^{\circ }$;
\item[$P_2)$] $v=v_{\mu _{q_1},q_{1}}^1$.\end{itemize} We observe
that pairs of this type exist also in $\cSN _B\times K^n$ (this
proves the observation written at the beginning of this section).
Since $v_{\mu _{q_1},q_{1}}^1\in \Delta _B ^{\circ }$, for these
pairs we have that $W_{X,v,1}= \langle \, \Delta _B ^{\circ }\,
\rangle $; hence by lemma \ref{nextt} we have that
$a(i,j,i',j',l')=0$ for all $i=1,\ldots ,\tilde i-1$,
$j=q_i-q_{i-1},\ldots ,1$ and $(i',j',l')\in H(i,j)$, that is
$\widehat {\Delta }_{B,X,v}= \{ W_{X,v,1}+v_{\mu _{q_i},j}^{l}\ |\
v_{\mu _{q_i},j}^{l}\in \Delta _B-\Delta _B^{\circ }\} $ by
Corollary \ref{next'}. If $v_{\mu _{q_i},j}^l$ and $v_{\mu
_{q_{i'}},j'}^{l'}$ are respectively the $h-$th and the $h'-$th
vector of $\Delta _B-\Delta _B^{\circ }$ then the entry of $p
_{X,v} \, \circ \, X_{v,1}\, \circ \, \big( p _{X,v}\big) ^{-1}\,
(e_h)\, $ with respect to $e_{h'}$ is the entry of $X\, \big(
v_{\mu _{q_i},j}^l\big) $ with respect to $v_{\mu
_{q_{i'}},j'}^{l'}$. By the definition of $R_B$, Proposition
\ref{order} and Corollary \ref{graph} we get that for all
$\widehat X\in \cSE _{\widehat B}$ there exists $(X,v)\in
\cSE^{\star  } _B$ with the properties $P_1)$ and $P_2)$ such that
$\pi _B ((X,v))=\widehat{X}$. \item[b)] Let $(X,v)\in \cQ _B$ and
let $(X',v')\in \cQ _B\, \cap \, \big( \cN _B \times K^n)$. By
Corollary \ref{t3} and Proposition \ref{j2} the endomorphisms
$X_{v,1}$ and $(X')_{v',1}$ are conjugated. Since $(X')_{v',1}$
commutes with $ J_{X',v'}$ we have that $\pi _B ((X',v'))\in \cN
_{\widehat B}$, hence $\pi _B ((X,v))\in \widetilde {\cSE }
_{\widehat
 B}$.\end{itemize} \hspace{4mm} $\square $\vspace{2mm}
\newline  Now we can prove Theorem \ref{eend}, which has been announced
in the introduction.\vspace{2mm} \newline {\bf Proof of Theorem
\ref{eend}}\vspace{2mm}\newline {\em $1^{st}$ step}:\
by Corollary \ref{t3} the maximum
partition which is associated to the elements of $\widetilde {\cSE
}_{\widehat B}\, $ ($\cSE _{\widehat B}$) is $(\widehat{\omega}_1,
\widehat {\omega}_2, \ldots ,\widehat{\omega}_{\widehat z})$.\vspace{2mm}\newline {\em $2^{nd}$ step}:\
 by Propositions \ref{f1} and by a) of Proposition \ref{end} there exists a non-empty
open subset of $\cSE _B$ such that if $X$ belongs to it the
partition of $X_{v,1}$ is greater or equal than $(\widehat{\omega}_1, \widehat {\omega}_2,
\ldots ,\widehat{\omega}_{\widehat z})$; hence by Proposition
\ref{j2} there exists a non-empty
open subset of $\cSE _B$ such that if $X$ belongs to it the partition of $X$ is greater or equal than
$(\omega _1,\widehat{\omega}_1, \widehat {\omega}_2,
\ldots ,\widehat{\omega}_{\widehat z})$.\vspace{2mm}\newline {\em $3^{rd}$ step}:\  by b) of Proposition \ref{end} the maximum partition
which is associated to the elements
of $\cSE _B$ is smaller or equal than $(\omega_1,\widehat{\omega}_1, \widehat{\omega}_2,
\ldots ,\widehat{\omega}_{\widehat z})$. \vspace{2mm}\newline Hence the maximum partition which is associated to the elements
of $\cSE _B$ is $(\omega_1,\widehat{\omega}_1, \widehat{\omega}_2,
\ldots ,\widehat{\omega}_{\widehat z})$, which by Corollary
\ref{t3} is equivalent to the claim. \hspace{4mm} $\square $\vspace{2mm}
\newline
{\bf Example 4} If $B=(15,13,5,4,3^2,2,1)$ we have $\tilde i=q_{\tilde i}=\mu _{q_{\tilde i}}=4$, $\omega_1=4+3\cdot 2+2\cdot 3=16$ and
$\widehat B=(15-2,13-2,5-2,2,1)=(13,11,3,2,1)$. Then $\widehat {\widehat B}=(11,3,2,1)$ and $\widehat {\widehat {\widehat B}}=(3,2,1)$. Since the maximum partition of the
elements of $\cN _{\widehat {\widehat {\widehat B}}}$ is $(5,1)$, we get that the maximum
partition of the elements of $\cN _B$ is $(16,13,11,5,1)$.\vspace{2mm}
\newline
Theorem \ref{eend} leads to an algorithm for the determination of
the maximum partition which is associated to elements of $\cN _B$
for any partition $B$.  \section{Other results on the orbits
intersecting ${\cN }_B$}  In \cite{BasI} one proved a relation
between $Q(B)$ and the Hilbert functions of the algebras $K[A,J]$
for $A\in \cN _B$. After this result Toma\v{z} Ko\v{s}ir and
Polona Oblak in \cite{Kos} (2008) proved that a generic pair of
commuting nilpotent matrices generates a Gorenstein algebra and
then obtained as a consequence the following result.
\begin{theorem}\label{KO} For any partition $B$ the partition
$Q(B)$ has decreasing parts differing by at least 2, hence the map
$Q$ is idempotent. \end{theorem} Theorem \ref{eend} leads to the
following proof of Theorem \ref{KO}. \vspace{2mm} \newline {\bf
Another proof of Theorem \ref{KO} } We can prove the claim by
induction on $n$, hence we can assume that $\widehat{\omega
}_{i}-\widehat{\omega }_{i+1}\geq 2$ for $i=1,\ldots ,\widehat
z-1$. By the definition of $\omega _1$ and $\widehat{\omega }_1$
we get that $\omega _1-\widehat {\omega }_1\geq 2$, hence by
Theorem \ref{eend} we get that $\omega _i-\omega _{i+1}\geq 2$ for
$i=1,\ldots ,z-1$. Then by Corollary \ref{corR} we get that
$Q(Q(B))=Q(B)$. \hspace{4mm} $\square $ \vspace{2mm}
\newline Polona Oblak recently published a paper on properties of
the nilpotent orbits which intersect $\cN _B$ for some special
types of $B$ (see \cite{Obll}), among which we cite the following
result.
\begin{theorem} \label{Obll} If $n>3$ then $\cN _B$ intersects all the nilpotent
orbits if and only if $J^2=0$.
\end{theorem}
The map $Q$ was investigated by D.I. Panyushev in the more general
context of Lie algebras. In fact if $\mathfrak{g}$ is a semisimple
Lie algebra over an algebraically closed field $K$ such that char
$K=0$ and $\cN (\mathfrak g)$ is the nilpotent cone of
$\mathfrak{g}$ then $\cN (\mathfrak g)$ is irreducible (see
\cite{Ko}).  Let $e\in \cN (\mathfrak{ g})$ and let
$\mathfrak{z}_{\mathfrak g}(e)$ be the centralizer of $e$. Let $\{
e,h,f\} $ be an $\mathfrak{sl}_2$-triple and let
$\mathfrak{g}={\displaystyle \bigoplus _{i\in \Z }\;
\mathfrak{g}(i)}$ be the corresponding $\Z-$ grading of
$\mathfrak{g}$.  Let $G$ be the adjoint group of $\mathfrak{g}$.
In \cite{Pan} D.I. Panyushev defined an element $e$ to be
self-large if $G\cdot e\ \cap \ (\mathfrak{z}_{\mathfrak g}(e)\,
\cap \, \cN (\mathfrak g))$ is open (dense) in
$\mathfrak{z}_{\mathfrak g}(e)\, \cap \, \cN (\mathfrak g)$ and he
 proved the following result.
\begin{theorem}\label{P} The element $e\in \cN $ is "self-large"
iff $\mathfrak{z}_{\mathfrak g}(e)\ \cap \ \mathfrak g (0)$ is
toral and $\mathfrak{z}_{\mathfrak g}(e)\ \cap \ \mathfrak g
(1)=\{ 0\} $. \vspace{4mm}\end{theorem}
 {\bf Acknowledgements} \ I am greatly indebted
to Polona Oblak for sharing the conjecture with me and to
Toma\v{z} Ko\v{s}ir for precious observations on the first version
of this paper.

{\small{Author's institution: Liceo Scientifico Annesso al
Convitto Nazionale, Assisi, Italy.}\newline {\small{Author's
address: Via dei Ciclamini 2B, 06126 Perugia Italy,
robasili@alice.it }
\end{document}